\documentclass[12pt,a4paper]{article}
\usepackage[a4paper,margin=1in]{geometry}
\usepackage{amsmath,amssymb,amsthm,mathtools}
\usepackage{microtype,booktabs,indentfirst}
\usepackage{xcolor}
\usepackage[colorlinks=true,linkcolor=blue,citecolor=blue,urlcolor=blue]{hyperref}

\numberwithin{equation}{section}
\newtheorem{theorem}{Theorem}[section]
\newtheorem{proposition}[theorem]{Proposition}
\newtheorem{lemma}[theorem]{Lemma}
\newtheorem{corollary}[theorem]{Corollary}
\newtheorem{conjecture}{Conjecture}

\theoremstyle{definition}
\newtheorem{definition}{Definition}[section]

\theoremstyle{definition}
\newtheorem{remark}[theorem]{Remark}

\begin{document}
	
\title{P\'olya's conjecture for hyperbolic balls and
	Berezin--Li--Yau type inequalities\thanks{{Research supported by the National Natural Science Foundation of China} (No. 12371110).}}
	
\author{
		Guowei Dai\thanks{School of Mathematical Sciences,
			Dalian University of Technology, Dalian 116024, P.R. China.
			Email: daiguowei@dlut.edu.cn.}
		\and
		Lan Zhao\thanks{Corresponding author. School of Mathematical Sciences,
			Dalian University of Technology, Dalian 116024, P.R. China.
			Email: zhaolan@mail.dlut.edu.cn.}}
	
\date{}
\maketitle

\begin{abstract}
We prove the strict shifted Dirichlet {P\'olya} inequality
\[
\lambda_k(B_R)>
\frac{(d-1)^2\kappa^2}{4}
+4\pi^2\left(
\frac{k}{\omega_d\operatorname{Vol}_\kappa(B_R)}
\right)^{2/d},
\qquad k\ge1,
\]
for every geodesic ball \(B_R\) in hyperbolic space of constant sectional curvature \(-\kappa^2\), in every dimension \(d\ge2\). As a consequence, we obtain the shifted {Berezin--Li--Yau} type inequality for hyperbolic balls, thereby providing a concrete family of domains satisfying both conjectures proposed by Laptev in 2023. 

For arbitrary nonempty open subsets \(\Omega\subset\mathbb{H}^d_\kappa\) of finite volume, we establish the weak {Berezin--Li--Yau} type inequality
\[
\sum_{j=1}^{m}\lambda_j(\Omega)
\ge
a_{d,\kappa}m
+\frac{d}{d+2}\,4\pi^2
\left(
\frac{m}{\omega_d\operatorname{Vol}_\kappa(\Omega)}
\right)^{2/d}m,
\qquad m\ge1,
\]
with the optimal leading coefficient, where \(a_{2,\kappa}=\kappa^2/4\) and \(a_{d,\kappa}=d(d-1)\kappa^2/6\) for \(d\ge3\). We also establish the global two-sided estimate for the renormalized eigenvalue sum on hyperbolic space, originally proposed by Strichartz in 1996 and valid for \(d\ge3\).
\end{abstract}


	
\section{Introduction}
Let $\Omega$ be a bounded nonempty open subset of Euclidean space $\mathbb{R}^d$, where $d\ge2$, and let $\left|\Omega\right|$ denote its volume. Consider the Dirichlet eigenvalue problem for the Laplacian
\[
\Delta := \sum_{j=1}^d \frac{\partial^2}{\partial x_j^2}
\]
in $\Omega$:
\begin{equation}
	\label{eq:dirichlet}
	\begin{cases}
		-\Delta u = \lambda u & \text{in } \Omega,\\
		u = 0 & \text{on } \partial\Omega.
	\end{cases}
\end{equation}
It is well known that the spectrum of \eqref{eq:dirichlet} is discrete and consists of isolated eigenvalues of finite multiplicity {accumulating at $+\infty$},
\[
0 < \lambda_1(\Omega) \le \lambda_2(\Omega) \le \cdots \le \lambda_n(\Omega) \le \dots,
\]
where each eigenvalue is repeated according to its multiplicity.

Weyl's law describes the asymptotic distribution of Dirichlet eigenvalues. For $\lambda\in\mathbb{R}$, define the eigenvalue counting function by
\begin{equation}\label{eq:intro-counting}
	N_\Omega\left(\lambda\right)=\#\left\{j\ge1:\lambda_j\left(\Omega\right)\le \lambda\right\},
\end{equation}
where $\#$ denotes cardinality, so each eigenvalue is counted according to its multiplicity. For smooth bounded domains, Weyl's law \cite{Chavel,Weyl} states that
\begin{equation}\label{eq:intro-weyl}
	N_\Omega\left(\lambda\right)
	=C_d\left|\Omega\right|\lambda^{d/2}+o\left(\lambda^{d/2}\right),
	\qquad \lambda\to+\infty.
\end{equation}
Here
\[
C_d=\frac{\omega_d}{\left(2\pi\right)^d}
\]
{is the so-called Weyl constant, and} $\omega_d$ is the volume of the Euclidean unit ball, given by
\[
\omega_d=\frac{\pi^{d/2}}{\Gamma\left(d/2+1\right)},
\]
where $\Gamma$ is the Gamma function.

P\'olya's conjecture asks whether the leading term in Weyl's law bounds the counting function {from above for every $\lambda>0$}. For a bounded nonempty open set $\Omega\subset\mathbb{R}^d$, the conjectured inequality is
\begin{equation}\label{eq:intro-polya-counting}
	N_\Omega\left(\lambda\right)\le C_d\left|\Omega\right|\lambda^{d/2},
	\qquad \lambda>0.
\end{equation}
Equivalently, it asserts that
\begin{equation}\label{eq:intro-polya-eigenvalues}
	\lambda_k\left(\Omega\right)\ge4\pi^2
	\left(\frac{k}{\omega_d\left|\Omega\right|}\right)^{2/d},
	\qquad k\ge1.
\end{equation}
Thus, while Weyl's law describes the asymptotic behavior as the eigenvalue tends to infinity, P\'olya's conjecture requires a bound for every eigenvalue.

P\'olya formulated the conjecture in \cite{Polya54} and proved the
Dirichlet case in \cite{Polya61} for Euclidean tiling domains $\Omega$: that is, domains such that $\mathbb{R}^d$ can be covered, up to a set of measure zero, {by a disjoint union of isometric copies of $\Omega$}. P\'olya also established the Neumann case under additional assumptions,
which were subsequently removed by Kellner \cite{Kellner}. Alongside direct proofs for particular classes of domains,
an important line of progress has been the study of averaged
spectral inequalities. Berezin's trace inequality \cite{Berezin}
and the Li--Yau inequality \cite{LiYau} provide such estimates
with sharp leading constants. In particular,
Li and Yau proved that the Dirichlet eigenvalues of a bounded
open set $\Omega\subset\mathbb{R}^d$ satisfy
\[
\sum_{j=1}^k\lambda_j\left(\Omega\right)
\ge \frac{d}{d+2}\,4\pi^2
\left(\omega_d\left|\Omega\right|\right)^{-2/d}
k^{1+2/d},
\qquad k\ge1.
\]
The leading constant is sharp by Weyl's law. However, the
resulting bound for an individual eigenvalue contains the factor
$d/(d+2)$, so it does not attain the bound asserted in {P\'olya's conjecture}.
Kr\"oger \cite{Kroger} obtained corresponding sharp estimates
for sums of Neumann eigenvalues.

Spectral asymptotics provide a complementary approach. For domains with piecewise smooth boundary, Weyl \cite{Weyl12} conjectured a two-term asymptotic formula whose second term reflects the boundary measure and the boundary condition.
Following H\"ormander's work on the spectral function
\cite{Hormander}, Ivrii \cite{Ivrii} established the second
term in Weyl's law under a dynamical assumption; see also
\cite{SV}. For a bounded domain with smooth boundary, if
periodic billiard trajectories form a set of measure zero
in phase space, Ivrii proved that the Dirichlet counting function satisfies
\begin{equation}\label{eq:intro-ivrii}
	N_\Omega\left(\lambda\right)
	=C_d\left|\Omega\right|\lambda^{d/2}
	-\frac{\omega_{d-1}}{4\left(2\pi\right)^{d-1}}
	\left|\partial\Omega\right|\lambda^{(d-1)/2}
	+o\left(\lambda^{(d-1)/2}\right),
	\qquad \lambda\to+\infty.
\end{equation}
Here $\left|\partial\Omega\right|$ denotes the surface measure
of the boundary. The negative boundary term implies
P\'olya's inequality for all sufficiently large eigenvalues, but does not control the remaining finite part
of the spectrum.

Beyond tiling domains, further progress has been made through product constructions. Laptev \cite{Laptev} proved that the Dirichlet
conjecture holds on $\Omega_1\times\Omega_2$ whenever it
holds on $\Omega_1\subset\mathbb{R}^{d_1}$ with $d_1\ge2$,
and $\Omega_2$ is an open set of finite volume.
More recently, He and Wang \cite{HeWang} established
the conjecture for sufficiently thin products
$(a\Omega_1)\times\Omega_2$, under the respective
hypotheses for Dirichlet and Neumann boundary conditions,
without assuming that either factor satisfies the
conjecture.

A major breakthrough in establishing P\'olya's conjecture for
specific non-tiling domains in Euclidean space was achieved by
Filonov, Levitin, Polterovich, and Sher \cite{FLPS}, who proved
the Dirichlet and Neumann conjectures for planar disks and
circular sectors, and the Dirichlet conjecture for Euclidean
balls in every dimension. Their method combines
estimates for Bessel functions with precise lattice
counting inequalities. In particular, the disk became
the first planar domain that does not tile the plane
for which the conjecture was verified. Their subsequent
work \cite{FLPSannuli} established the Dirichlet conjecture
for all planar concentric annuli. The recent preprint of Tang and Wang \cite{TW} proves a strict Dirichlet inequality for bounded simply connected planar domains. The conjecture for arbitrary bounded Euclidean
domains remains open.

In this paper, we consider the complete simply connected hyperbolic space $\mathbb{H}^d_\kappa$ with sectional curvature $-\kappa^2$, where $d\ge 2$ and $\kappa>0$. Let $\operatorname{Vol}_\kappa$ denote its Riemannian volume. Its nonnegative Laplace--Beltrami operator has spectral bottom
\[
\rho_{d,\kappa}^2=\frac{\left(d-1\right)^2\kappa^2}{4};
\]
see McKean \cite{McKean}, Helgason \cite{Helgason}, and
\cite[Proposition~2.1]{ILW}. In his 2023 lecture
\textit{Spectral inequalities for Schr\"odinger operators on manifolds with constant curvature} \cite{Laptev2023},
Laptev stated two conjectures for the
Dirichlet spectrum after subtracting the spectral bottom.
We state them below in our notation, with curvature scaled from $-1$
to $-\kappa^2$. 

\begin{conjecture}[Shifted P\'olya inequality; Laptev]\label{conj:laptev-polya}
Let $\Omega\subset\mathbb H^d_\kappa$ be a nonempty open set of finite
volume, where $d\ge2$. For every $s>0$,
\begin{equation}\label{eq:hyperbolic-target}
N_\Omega(\rho_{d,\kappa}^2+s^2)
\le C_d\operatorname{Vol}_\kappa(\Omega)s^d.
\end{equation}
\end{conjecture}

\begin{conjecture}[Shifted Berezin--Li--Yau inequality; Laptev]\label{conj:laptev-bly}
Under the same assumptions, for every $s>0$,
\begin{equation}\label{eq:laptev-bly-target}
R_\Omega(\rho_{d,\kappa}^2+s^2)
\le L^{\mathrm{cl}}_{1,d}\operatorname{Vol}_\kappa(\Omega)s^{d+2},
\end{equation}
where 
$$
R_\Omega(\lambda)=\sum_{j=1}^\infty(\lambda-\lambda_j(\Omega))_+
$$ 
and
\[
x_+=\max\{x,0\},
\qquad
L^{\mathrm{cl}}_{1,d}
=\frac{1}{(4\pi)^{d/2}\Gamma(d/2+2)}=\frac{2}{d+2}C_d.
\]
Equivalently, for every $m\in\mathbb N$,
\[
\sum_{j=1}^m\lambda_j(\Omega)
\ge\rho_{d,\kappa}^2m+
\frac{d}{d+2}\,4\pi^2
\left(\frac{m}{\omega_d\operatorname{Vol}_\kappa(\Omega)}\right)^{2/d}m.
\]
\end{conjecture}

Ilyin, Laptev, and Weinmann studied shifted spectral estimates in \cite{ILW} and revived the first conjecture. Their work makes a substantial contribution to the theory of spectral inequalities in hyperbolic space and provides an important foundation for the problems considered here. The first conjecture implies the second by integration of the counting function. The converse does not follow from this averaging argument.
One contribution of the present work is to verify both conjectures
for every hyperbolic geodesic ball, in every dimension $d\ge2$.
Our counting estimate is strict and includes an explicit positive
remainder. We also obtain a version of Conjecture~\ref{conj:laptev-bly}
for arbitrary open sets of finite volume in every dimension $d\ge2$,
with the same leading coefficient. The full conjectured shift is
obtained in dimensions two and three, and a weaker curvature shift
is obtained in dimensions $d\ge4$.

Our first main result establishes the strict shifted inequality for every geodesic ball and estimates the remainder from the Weyl term. For $z>0$, define
\begin{equation}\label{eq:epsilon-definition}
	\varepsilon\left(z\right)
	=\frac1\pi\int_0^\infty
	\frac{\log\left|\dfrac{z+t}{z-t}\right|}{e^{2\pi t}-1}\,\mathrm{d} t.
\end{equation}
The singularity at $t=z$ is logarithmic and integrable. We will prove that $0<\varepsilon\left(z\right)<1/4$. For $\nu \ge0$, we define
\begin{equation}\label{eq:general-action}
	G_{\kappa,R,s}\left(\nu \right)
	=\frac1\pi\int_0^R
	\sqrt{\left(s^2-\kappa^2\nu^2\operatorname{csch}^2\left(\kappa r\right)\right)_+}\,\mathrm{d} r{.}
\end{equation}
The positive remainder is
\begin{equation}\label{eq:ball-deficit}
	\mathcal D_{d,\kappa}\left(R,s\right)
	=\frac2{\left(d-2\right)!}\int_0^{\nu_{\max,\kappa,R,s}}
	{\nu^{d-2}\min\left\{G_{\kappa,R,s}\left(\nu \right),\varepsilon\left(s/\kappa\right)\right\}\,\mathrm{d} \nu,}
\end{equation}
where
\[
\nu_{\max,\kappa,R,s}=\frac{s}{\kappa}\sinh\left(\kappa R\right).
\]
These quantities depend only on $d$, $\kappa$, $R$, and $s$, and do not involve any unknown eigenvalues.

\begin{theorem}\label{thm:balls}
	For every geodesic ball $B_R\subset\mathbb{H}^d_\kappa$, every $R>0$ and every $s>0$, the remainder in \eqref{eq:ball-deficit} satisfies
	\[
	0<\mathcal D_{d,\kappa}\left(R,s\right)
	\le C_d\operatorname{Vol}_\kappa\left(B_R\right)s^d.
	\]
	The counting function satisfies
	\begin{equation}\label{eq:ball-count-refined}
		N_{B_R}\left(\rho_{d,\kappa}^2+s^2\right)
		\le C_d\operatorname{Vol}_\kappa\left(B_R\right)s^d-\mathcal D_{d,\kappa}\left(R,s\right).
	\end{equation}
	The inequality in \eqref{eq:ball-count-refined} is strict whenever its right side is positive. In particular,
	\begin{equation}\label{eq:ball-count}
		N_{B_R}\left(\rho_{d,\kappa}^2+s^2\right)
		<C_d\operatorname{Vol}_\kappa\left(B_R\right)s^d.
	\end{equation}
	Consequently, for all $k\in\mathbb{N}$,
	\begin{equation}\label{eq:ball-eigen}
		\lambda_k\left(B_R\right)>
		\rho_{d,\kappa}^2+4\pi^2
		\left(\frac{k}{\omega_d\operatorname{Vol}_\kappa\left(B_R\right)}\right)^{2/d}.
	\end{equation}
	The leading coefficient and the constant spectral shift in \eqref{eq:ball-eigen} are individually optimal uniformly over all radii and spectral indices.
\end{theorem}

Theorem~\ref{thm:balls} proves Conjecture~\ref{conj:laptev-polya}
for balls. Its integrated consequence, Corollary~\ref{cor:ball-bly},
proves Conjecture~\ref{conj:laptev-bly} for the same class of domains.
Thus geodesic balls give a concrete family satisfying both conjectures,
with the full spectral shift and the sharp Euclidean coefficients.
Corollaries~\ref{cor:sectors}--\ref{cor:tiles} give the corresponding strict inequalities for planar sectors, disjoint unions of balls with finite total volume, and congruent subdivisions of balls.

We now outline the main ideas behind the proof of
Theorem~\ref{thm:balls}. The main difficulty is to obtain a counting estimate that preserves both the sharp volume coefficient and the spectral shift
$\rho_{d,\kappa}^2$. Separation of variables reduces the problem
to a family of radial equations with potentials involving
$1/\sinh^2(\kappa r)$. To recover the full counting function,
one must control the zeros of the radial solutions for every
angular order $\ell\in\mathbb{N}_0$ and then sum the resulting bounds with the
corresponding multiplicities. Errors in these bounds may accumulate
under summation and prevent the desired sharp estimate.
We address this difficulty by deriving a radial phase bound
that is valid for every spectral parameter, ball radius, and
angular order and is sufficiently precise to retain a positive
remainder after summation. This uniform
control allows us to apply the lattice inequalities
of \cite{FLPS} without losing the sharp volume
coefficient, while the correction yields an explicit
remainder and makes the resulting P\'olya inequality
strict.

Theorem~\ref{thm:balls} establishes the shifted P\'olya
conjecture for every hyperbolic geodesic ball in every dimension. The rotational symmetry of the ball is essential to
our proof because it allows separation of variables and reduces the counting problem to radial phase estimates and lattice counting.
For a general domain, the Dirichlet boundary condition is not compatible with this separation, so the radial equations and the associated counting argument are no longer directly available.
Therefore, our second main result addresses arbitrary nonempty
open subsets of $\mathbb H^d_\kappa$ of finite volume through
a different approach, combining the spectral resolution with localized trace estimates. We set
\begin{equation}\label{eq:averaged-shift}
a_{d,\kappa}=
\begin{cases}
\kappa^2/4,&d=2,\\[2pt]
d(d-1)\kappa^2/6,&d\ge3.
\end{cases}
\end{equation}

\begin{theorem}\label{thm:averaged}
Let $d\ge2$, and let $\Omega\subset\mathbb H^d_\kappa$ be a nonempty
open set with $V=\operatorname{Vol}_\kappa(\Omega)<\infty$. For every $\lambda\in\mathbb R$, we have
\begin{equation}\label{eq:general-riesz}
R_\Omega(\lambda)
\le L^{\mathrm{cl}}_{1,d}V
(\lambda-a_{d,\kappa})_+^{1+d/2}.
\end{equation}
Equivalently, for every $m\in\mathbb N$,
\begin{equation}\label{eq:sum-inequality}
\sum_{j=1}^m\lambda_j(\Omega)
\ge a_{d,\kappa}m+
\frac{d}{d+2}\,4\pi^2
\left(\frac{m}{\omega_d V}\right)^{2/d}m.
\end{equation}
In particular, when $d=2$, these bounds become
\begin{equation}\label{eq:planar-riesz}
R_\Omega(\lambda)
\le\frac{V}{8\pi}\left(\lambda-\frac{\kappa^2}{4}\right)_+^2
\end{equation}
and
\begin{equation}\label{eq:planar-sum}
\sum_{j=1}^m\lambda_j(\Omega)
\ge\frac{\kappa^2}{4}m+\frac{2\pi}{V}m^2.
\end{equation}
The leading coefficient in the eigenvalue sum bound cannot be increased.
In dimension two the constant shift $\kappa^2/4$ is also optimal.
\end{theorem}
\begin{remark}
We note that the three-dimensional eigenvalue sum
inequality in Theorem~\ref{thm:averaged} is already known.
For bounded open sets, it follows from
\cite[Theorem~4.2 and equation~(5.15)]{Strichartz}
after scaling the curvature to $-\kappa^2$.	
\end{remark}

Conjecture~\ref{conj:laptev-bly} asks for
\eqref{eq:general-riesz} and \eqref{eq:sum-inequality} with $a_{d,\kappa}$ replaced by
$\rho_{d,\kappa}^2$. The two shifted bounds coincide when $d=2,3$, whereas
\[
\rho_{d,\kappa}^2-a_{d,\kappa}
=\frac{(d-1)(d-3)\kappa^2}{12}>0
\qquad d\ge4.
\]
This gap reflects a limitation of the trace approach. After
subtracting the spectral bottom, the spectral density in
dimensions $d\ge4$ exceeds its Euclidean counterpart at large
spectral parameters, so the trace estimate does not yield
the full conjectured shift.
Thus Theorem~\ref{thm:averaged} establishes the conjectured averaged
inequality in dimensions two and three and advances Conjecture~\ref{conj:laptev-bly}
in higher dimensions by retaining the sharp leading coefficient and
a weaker positive curvature shift. These averaged bounds are also weaker than the individual eigenvalue estimates in Theorem~\ref{thm:balls} and do not by themselves establish the P\'olya conjecture for general domains.
Nevertheless, they provide quantitative progress towards
the conjectured inequalities for general hyperbolic domains. For balls, Corollary~\ref{cor:ball-bly} gives the full shift in every
dimension $d\ge2$.
In dimension two, \eqref{eq:planar-sum} extends the Euclidean
Li--Yau inequality by the optimal curvature correction
$\kappa^2m/4$.

Our estimates for eigenvalue sums on general hyperbolic domains
also connect with Strichartz's work
on eigenvalue sums for domains in homogeneous spaces \cite{Strichartz}.
Strichartz used the spectral resolution of the ambient Laplacian to construct
a renormalized eigenvalue sum $\Sigma_{d,\kappa}(x)$. The precise definition is
given in Section~\ref{subsec:strichartz}. Theorem~5.1 in \cite{Strichartz} identifies
$d(d-1)/6$ as the linear coefficient as $x\to\infty$
at curvature $-1$, and $(d-1)^2/4$ as the linear coefficient as $x\rightarrow0$.
{These asymptotic formulas do not provide bounds valid for
every $x>0$.} In \cite[equation~(5.24)]{Strichartz}, Strichartz proposed a global two-sided
bound and conjectured these two coefficients as the optimal constants
when $d\ge3$. Strichartz left this global estimate unproved in \cite{Strichartz}, and, to the best
of our knowledge, no proof has appeared in the subsequent literature. {Proposition~\ref{prop:strichartz} establishes this bound.} Specifically, for every $d\ge3$ and every $x>0$,
\begin{equation}\label{eq:strichartz-intro}
a_{d,\kappa}x
\le\Sigma_{d,\kappa}(x)
-\frac{d}{d+2}\,4\pi^2\omega_d^{-2/d}x^{1+2/d}
\le\rho_{d,\kappa}^2x.
\end{equation}
Both linear coefficients are optimal. At curvature $-1$, this proves
the global estimate proposed in \cite[equation~(5.24)]{Strichartz}
with the conjectured constants.
The lower coefficient $a_{d,\kappa}$ also provides the shift in
Theorem~\ref{thm:averaged}, while the upper coefficient
$\rho_{d,\kappa}^2$ is the shift in Conjecture~\ref{conj:laptev-bly}.
These two claims must be distinguished: an upper bound on the renormalized sum does not yield the conjectured lower bound for Dirichlet eigenvalue sums. Consequently, optimality of $a_{d,\kappa}$ in \eqref{eq:strichartz-intro} does not imply its optimality in \eqref{eq:sum-inequality}. Thus, Conjecture~\ref{conj:laptev-bly} remains open for general
open sets of finite volume in dimensions $d\ge4$.

The paper is organized as follows.
Section~\ref{sec:prelim} collects the geometric preliminaries,
lattice inequalities, radial equations, special function identities,
and localized trace formulas. Section~\ref{sec:balls} proves Theorem~\ref{thm:balls} and its
corollaries, including the shifted Berezin--Li--Yau inequality for
balls. Section~\ref{sec:general} proves Theorem~\ref{thm:averaged}
and Proposition~\ref{prop:strichartz}.

\section{Notation and preliminaries}\label{sec:prelim}
We first introduce the basic notation and recall the geometric and spectral facts used throughout the paper. 

\subsection{Geometric and spectral preliminaries}

Let $\left(M,g\right)$ be a smooth Riemannian manifold and let
$\Omega\subset M$ be a nonempty open set. The Riemannian volume element is $\,\mathrm{d} V_g$ and the metric gradient is $\nabla_g$. The space $H^1_0\left(\Omega\right)$ is the closure of $C_c^\infty\left(\Omega\right)$ {in the norm whose square is
}\[
\left\|u\right\|_{H^1(\Omega)}^2=\int_\Omega\left(\left|u\right|^2+\left|\nabla_g u\right|_g^2\right)\,\mathrm{d} V_g.
\]
{Here $C_c^\infty(\Omega)$ consists of smooth functions compactly supported in $\Omega$.} The $L^2(\Omega)$ inner product and norm are defined by
\[
\langle u,v\rangle_{L^2(\Omega)}
=\int_\Omega u\overline{v}\,\mathrm{d}V_g,
\qquad
\|u\|_{L^2(\Omega)}^2
=\int_\Omega |u|^2\,\mathrm{d}V_g.
\] 
The Dirichlet Laplacian $-\Delta_{\Omega}$ is the nonnegative self-adjoint operator on $L^2(\Omega)$ associated with the closed form
\[
q_\Omega\left[u\right]=\int_\Omega\left|\nabla_g u\right|_g^2\,\mathrm{d} V_g,
\qquad u\in H^1_0\left(\Omega\right).
\]

For hyperbolic space, $g$ denotes the hyperbolic metric and
$\mathrm dV_\kappa$ its volume element. In dimension two we also write
$\mathrm dA_\kappa=\mathrm dV_\kappa$.
On the whole space let $H=-\Delta_{\mathbb H^d_\kappa}$, whose form is
\[
q_H[v]=\int_{\mathbb H^d_\kappa}|\nabla_gv|_g^2\,\mathrm dV_\kappa,
\qquad D(q_H)=H^1(\mathbb H^d_\kappa).
\]
For $u\in H^1_0(\Omega)$, its zero extension $\widetilde u$ belongs
to $H^1(\mathbb H^d_\kappa)$ and satisfies
$q_H[\widetilde u]=q_\Omega[u]$.
For every nonempty hyperbolic open set of finite volume, the
Dirichlet spectrum is discrete; see \cite{ILW}. All min--max and counting statements
below are applied to such sets, or to relatively compact open sets.

The min--max principle gives
\[
\lambda_k\left(\Omega\right)
=\inf_{\substack{L\subset H^1_0\left(\Omega\right)\\\dim L=k}}
\sup_{0\ne u\in L}\frac{q_\Omega\left[u\right]}{\left\|u\right\|_{L^2(\Omega)}^2}.
\]
We use its standard consequences, including domain monotonicity and Dirichlet bracketing \cite{Davies}. In particular, if $D_1,\ldots,D_q$ are disjoint open subsets of $\Omega$, zero extension gives
\begin{equation}\label{eq:bracketing}
	\sum_{i=1}^qN_{D_i}\left(\lambda\right)\le N_\Omega\left(\lambda\right).
\end{equation}
When $\Omega$ itself is their disjoint union, the operator is the orthogonal direct sum and the counts add exactly.

The spectrum on the whole space \(\mathbb{H}^d_\kappa\) is
\[
\left[\rho_{d,\kappa}^2,\infty\right).
\]
In particular, $\rho_{d,\kappa}^2$ is not an $L^2$ eigenvalue \cite{ILW} and
\begin{equation}\label{spectrum-bottom}
\lambda_j\left(\Omega\right)>\rho_{d,\kappa}^2.
\end{equation}

\begin{definition}\label{def:count-riesz}
For $\lambda\in\mathbb{R}$, the first Riesz mean is defined by
\[
R_\Omega(\lambda)
=\sum_{j\ge1}\left(\lambda-\lambda_j(\Omega)\right)_+.
\]
\end{definition}
It is easy to see that only eigenvalues strictly below $\lambda$ contribute to the sum, so \(R_\Omega(\lambda)\) is finite for fixed $\lambda$. We use $\mathbb{N}=\left\{1,2,\ldots\right\}$ and $\mathbb{N}_0=\left\{0,1,2,\ldots\right\}$. {The notation} $\lfloor x\rfloor=\max\{k\in\mathbb Z:k\le x\}$
denotes the integer part of $x\in\mathbb R$, and
$\lceil x\rceil=\min\{k\in\mathbb Z:k\ge x\}$
denotes its ceiling.

The metric of \(\mathbb{H}^d_\kappa\) in polar coordinates is
\begin{equation}\label{eq:polar}
	\,\mathrm{d} r^2+\left(\frac{\sinh\left(\kappa r\right)}{\kappa}\right)^2\,\mathrm{d}\omega^2,
\end{equation}
where $r$ is geodesic distance from the center and $\,\mathrm{d}\omega^2$ {is the standard metric on the unit sphere $\mathbb{S}^{d-1}$.} A geodesic ball is denoted by $B_R$, with $R$ its radius. Its volume is
\begin{equation}\label{eq:ball-volume}
	\operatorname{Vol}_\kappa\left(B_R\right)
	=\left|\mathbb{S}^{d-1}\right|\int_0^R\left(\frac{\sinh\left(\kappa r\right)}{\kappa}\right)^{d-1}\,\mathrm{d} r.
\end{equation}
Here the area of the unit sphere is
\[
\left|\mathbb{S}^{d-1}\right|=\frac{2\pi^{d/2}}{\Gamma\left(d/2\right)}.
\]

\subsection{Lattice inequalities}
We recall the two-dimensional lattice theorem of
Filonov, Levitin, Polterovich, and Sher \cite{FLPS}
and then state its consequence in higher dimensions. 

\begin{theorem}[See {\cite[Theorem~5.1]{FLPS}}]\label{thm:lattice-two}
	Let $\nu_{\max}>0$. Suppose that $g:\left[0,\nu_{\max}\right]\to\left[0,\infty\right)$ is nonincreasing, convex, {and $1/2$-Lipschitz, with $g\left(\nu_{\max}\right)=0$}. Then
	\begin{equation}\label{eq:lattice-two}
		\left\lfloor g\left(0\right)+\frac14\right\rfloor
		+2\sum_{j=1}^{\left\lfloor \nu_{\max}\right\rfloor}\left\lfloor g\left(j\right)+\frac14\right\rfloor
		\le2\int_0^{\nu_{\max}} g\left(\nu \right)\,\mathrm{d} \nu .
	\end{equation}
	Equality holds if and only if $g$ is identically zero.
\end{theorem}

For $d\ge3$, the dimension of the space of spherical harmonics of degree $\ell$ on $\mathbb{S}^{d-1}$ {(see \cite{Chavel})} is
\[
{m_{d,\ell}}=\binom{\ell+d-1}{d-1}-\binom{\ell+d-3}{d-1},
\qquad \ell\in\mathbb{N}_0.
\]
We use the convention that $\binom{m}{k}=0$
for integers $m<k$ with $k\ge0$. For $d=2$, the corresponding multiplicities are
${m_{2,0}}=1$ and ${m_{2,\ell}}=2$ for $\ell\ge1$. Thus Theorem~\ref{thm:lattice-two} may be regarded as the \(d=2\) case of the theorem below.

\begin{theorem}\label{thm:lattice-high}
	Let $d\ge3$ and $\nu_{\max}>0$. Suppose that $g:\left[0,\nu_{\max}\right]\to\left[0,\infty\right)$ is nonincreasing, convex, {and $1/2$-Lipschitz, with $g\left(\nu_{\max}\right)=0$ and $g\left(\nu\right)>0$ for $0\le\nu<\nu_{\max}$.} We extend $g$ by zero for \(\nu >\nu_{\max}\). Then
	\begin{equation}\label{eq:lattice-weighted}
		\sum_{\ell\ge0}{m_{d,\ell}}
		\left\lfloor g\left(\ell+\frac{d-2}{2}\right)+\frac14\right\rfloor
		<\frac2{\left(d-2\right)!}\int_0^{\nu_{\max}} \nu ^{d-2}g\left(\nu \right)\,\mathrm{d} \nu .
	\end{equation}
\end{theorem}
\begin{proof}
	Set 
	$$
	h=\frac{d-2}{2},\quad  c_n=\binom{n+d-3}{d-3},\quad  q_\ell=\left\lfloor g(\ell+h)+\frac14\right\rfloor.
	$$
	{If $\nu_{\max}\le h$, \eqref{eq:lattice-weighted} is trivial, so we suppose $\nu_{\max}>h$.} By \cite[Theorem~7.1]{FLPS}, we have
	$$
	{m_{d,\ell}}=c_\ell+2\sum_{n=0}^{\ell-1}c_n.
	$$
	For each integer $n\ge 0$ satisfying $n+h<\nu_{\max}$, {we apply Theorem~\ref{thm:lattice-two} to the function $t\mapsto g(t+n+h)$ on $[0,\nu_{\max}-n-h]$ to obtain}
	\[
	q_{n}+2\sum_{j\ge 1}q_{n+j}<2\int_{n+h}^{\nu_{\max}}g(\nu)\,\mathrm{d}\nu.
	\]
	Here the left-hand side corresponds to the evaluations of the function $t\mapsto g(t+n+h)$ at $t=0,1,2,\dots$, and the right-hand side follows from the change of variables $\nu = t+n+h$. Thus, we obtain
	\begin{align*}
		\sum_{\ell\ge0}{m_{d,\ell}}q_\ell
		&=\sum_{n+h<\nu_{\max}}c_n\left(q_n+2\sum_{j\ge1}q_{n+j}\right)\\
		&<2\sum_{n+h<\nu_{\max}}c_n\int_{n+h}^{\nu_{\max}} g(\nu)\,\mathrm{d}\nu\\
		&=2\int_0^{\nu_{\max}} f_d(\nu)g(\nu)\,\mathrm{d}\nu,
	\end{align*}
	{where
	\[
	f_d(\nu)=\sum_{\substack{n\in\mathbb N_0\\n+h\le\nu}}c_n
	=\begin{cases}
	0, & 0\le\nu<h,\\
	\displaystyle\binom{\lfloor\nu-h\rfloor+d-2}{d-2}, & \nu\ge h,
	\end{cases}
	\]
	is the function in \cite[Equation~(7.4)]{FLPS}.} By \cite[Equation~(7.6)]{FLPS}, we have
	\[
	\int_0^z f_d(\nu)\,\mathrm{d}\nu
	\le \frac{z^{d-1}}{(d-1)!},\qquad z\ge0.
	\]
	Since $g$ is nonnegative and nonincreasing,
	\cite[Lemma~7.2]{FLPS} gives
	\[
	\int_0^{\nu_{\max}} f_d(\nu)g(\nu)\,\mathrm{d}\nu
	\le \frac{1}{(d-2)!}\int_0^{\nu_{\max}} \nu^{d-2}g(\nu)\,\mathrm{d}\nu.
	\]
	Combining these estimates proves \eqref{eq:lattice-weighted}.
\end{proof}

\subsection{The radial equation and zero counting}
\label{subsec:singular-oscillation}

In this subsection, we derive the radial equation for hyperbolic balls, characterize the solution regular at the origin, and recall the corresponding counting theorem.

We take $\kappa=1$ and set $\rho=\left(d-1\right)/2$. A spherical harmonic of degree $\ell\in\mathbb{N}_0$ is an eigenfunction $Y_\ell$ of the operator $-\Delta_{\mathbb{S}^{d-1}}$ with eigenvalue $\ell\left(\ell+d-2\right)$. The Laplacian corresponding to \eqref{eq:polar} is
\[
\Delta=\partial_r^2+\left(d-1\right)\coth r\,\partial_r
+\sinh^{-2}r\,\Delta_{\mathbb{S}^{d-1}}.
\]
We set $u\left(r,\omega\right)=\left(\sinh r\right)^{-\rho}w\left(r\right)Y_\ell\left(\omega\right)$. Substituting $u\left(r,\omega\right)$ into the eigenvalue equation and setting $\lambda=\rho^2+s^2$ and $\nu =\ell+\left(d-2\right)/2$, {we obtain}
\begin{equation}\label{eq:radial}
	-w''+\frac{\nu^2-1/4}{\sinh^2r}w=s^2w \qquad 0<r<R
\end{equation}
and $w\left(R\right)=0$.

{To select the appropriate solution, we want to determine its leading behavior near $r=0$.} Since $\sinh r\sim r$ as $r\rightarrow0$, {substituting $w(r)= r^\beta$ into
\eqref{eq:radial} and letting $r\rightarrow0$, we obtain}
\[
-\beta(\beta-1)+\nu^2-\frac14=0.
\]
For $\nu>0$, the roots of this equation are $1/2+\nu$ and $1/2-\nu$. When $\nu=0$, the root is repeated and the two independent
solutions have leading terms $r^{1/2}$ and
$r^{1/2}\log r$, respectively.	The solution we seek must satisfy a suitable regularity condition at the origin. For the radial equation considered here, {regularity at the origin is expressed by the following Friedrichs condition;}
see \cite{GPS}.
\begin{definition}\label{def:friedrichs-origin}
A nonzero solution $w$ of the radial equation is said to
satisfy the Friedrichs condition at the origin if
\begin{equation}\label{equ:friedrichs-origin}
w(r)=c\,r^{1/2+\nu}\bigl(1+O(r^2)\bigr)
\qquad \text{as }r\rightarrow0,	
\end{equation}
where $c\ne0$ is a constant.	
\end{definition}

{To relate this condition to regularity of the original eigenfunction, note that} for $\nu=\ell+(d-2)/2$ and $\rho=(d-1)/2$,
the radial part satisfies
\[
(\sinh r)^{-\rho}w(r)
=c\,r^\ell\left(1+O(r^2)\right),
\]
which is compatible with regularity at $r=0$.
{For $\nu>0$, the other case of the radial part has leading behavior
$r^{-(\ell+d-2)}$ and is excluded. When $\nu=0$, the excluded solution
$w(r)\sim\sqrt r\log r$ corresponds to logarithmic behavior of the original radial part. See \cite{NZ,Zettl}.}

For any real $\nu\ge0$, we define
\begin{equation}\label{eq:pre-sl-expression}
\tau_\nu y=-y''+\frac{\nu^2-1/4}{\sinh^2r}y,
\qquad 0<r<R.
\end{equation}
Let \(L_{\nu,R}\) denote its Friedrichs realization in \(L^2(0,R)\), subject to \(y(R)=0\) and the origin condition described above. For ${\mu}\in\mathbb R$, let $y_\nu(r,{\mu})$ solve
$\tau_\nu y={\mu} y$ normalized by
\begin{equation}\label{eq:pre-principal}
	\lim_{r\rightarrow0}\frac{y_\nu(r,{\mu})}{r^{1/2+\nu}}=1.
\end{equation}
{We have the following result.}
\begin{lemma}
	\label{lem:pre-oscillation}
	The operator $L_{\nu,R}$ has simple discrete eigenvalues
	$\mu_{\nu,1}<\mu_{\nu,2}<\cdots$ tending to $+\infty$.
	{We have the counting identities}
	\begin{equation}\label{eq:pre-oscillation-strict}
		\#\{j:\mu_{\nu,j}<{\mu}\}
		=\#\{r\in(0,R):y_\nu(r,{\mu})=0\},
	\end{equation}
	and
	\begin{equation}\label{eq:pre-oscillation-inclusive}
		\#\{j:\mu_{\nu,j}\le {\mu}\}
		=\#\{r\in(0,R]:y_\nu(r,{\mu})=0\}.
	\end{equation}
	In particular, the $j$th eigenfunction has exactly $j-1$
	zeros in $(0,R)$.
\end{lemma}

{For completeness, we outline the proof.} The Friedrichs realization of the following operator
$$
-\frac{\mathrm{d}^2}{\mathrm{d}r^2}+\frac{\nu^2-\frac14}{r^2}
$$
on $(0,R)$, with the Dirichlet condition at $R$, has
compact resolvent. Since
$(\nu^2-\frac14)(\sinh^{-2}r-r^{-2})$ is bounded,
$L_{\nu,R}$ also has compact resolvent by stability
under bounded perturbations. Consequently, its spectrum
is purely discrete, with eigenvalues tending to $+\infty$. Simplicity follows from the fact that solutions satisfying the Friedrichs condition at the origin are unique up to a constant multiple. The zero count follows from \cite[Theorem~2.1,
Corollaries~2.3 and 2.4]{GST}. The formula with $\mu_{\nu,j}\le {\mu}$ includes the zero at $R$ if and only if ${\mu}$ is an eigenvalue.

\begin{remark}
{Suppose that $y_\nu(r,\mu)=A(r)\sin\beta(r)$, where $A(r)>0$ for $r\in(0,R]$, and that $\beta$ is continuous and strictly increasing on $[0,R]$, with $\beta(0)=0$.} Then \eqref{eq:pre-oscillation-inclusive} gives
\begin{equation}\label{eq:pre-oscillation-phase}
	\#\{j:\mu_{\nu,j}\le {\mu}\}
	=\left\lfloor\frac{\beta(R)}{\pi}\right\rfloor.
\end{equation}
{For the ball problem, the radial parameter is $\mu=s^2=\lambda-\rho^2$. We will use this relation later.}	
\end{remark}

\subsection{Special functions}\label{subsec:special-oscillation}

In what follows, we shall use several basic properties of the Gamma,
{Beta and Gauss hypergeometric functions.} For the reader's convenience,
we recall their definitions and the identities needed in
Section~\ref{sec:balls}. For proofs, we refer to
\cite[Chapters~5 and 15]{DLMF}.

\paragraph[Gamma and beta functions.]{{Gamma and Beta functions.}}
For $\operatorname{Re}z>0$, the Gamma function is defined by
\begin{equation}\label{eq:pre-gamma}
\Gamma(z)=\int_0^\infty t^{z-1}e^{-t}\,\mathrm{d}t.
\end{equation}
Its recurrence is
\begin{equation}\label{eq:pre-gamma-recurrence}
\Gamma(z+1)=z\Gamma(z).
\end{equation}
In particular, $\Gamma(n+1)=n!$ for $n\in\mathbb N_0$ and
$\Gamma(1/2)=\sqrt\pi$. We also use the reflection formula
\begin{equation}\label{eq:pre-reflection}
	\Gamma(z)\Gamma(1-z)
	=\frac{\pi}{\sin(\pi z)},
	\qquad z\in\mathbb C\setminus\mathbb Z.
\end{equation}
Combining the recurrence relation and the reflection formula, we obtain, for $a>0$,
\begin{equation}\label{eq:pre-gamma-imaginary-modulus}
	|\Gamma(\mathrm ia)|^2
	=\frac{\pi}{a\sinh(\pi a)},
\qquad
	|\Gamma(1/2+\mathrm ia)|^2
	=\frac{\pi}{\cosh(\pi a)}.
\end{equation}
The duplication formula is
\begin{equation}\label{eq:pre-duplication}
\Gamma(z)\Gamma(z+1/2)=2^{1-2z}\sqrt\pi\,\Gamma(2z),
\qquad \operatorname{Re}z>0.
\end{equation}
For $\operatorname{Re}x>0$ and $\operatorname{Re}y>0$, the beta function is
\begin{equation}\label{eq:pre-beta}
\mathrm B(x,y)=\int_0^1t^{x-1}(1-t)^{y-1}\,\mathrm{d}t.
\end{equation}
Its relation to the Gamma function is
\begin{equation}\label{eq:pre-beta-gamma}
\mathrm B(x,y)=\frac{\Gamma(x)\Gamma(y)}{\Gamma(x+y)}.
\end{equation}

Throughout the paper, $\operatorname{Log}\Gamma$ in the right half-plane
is the analytic branch that is real on the positive real axis. $\operatorname{Log}z$ is the principal logarithm. Thus
\begin{equation}\label{eq:pre-log-recurrence}
\operatorname{Log}\Gamma(z+1)
=\operatorname{Log}\Gamma(z)+\operatorname{Log}z,
\qquad \operatorname{Re}z>0.
\end{equation}
Stirling's formula on this branch is
\begin{equation}\label{eq:pre-stirling}
\operatorname{Log}\Gamma(z)
=\left(z-\frac12\right)\operatorname{Log}z-z
+\frac12\log(2\pi)+O(|z|^{-1}),
\end{equation}
as $|z|\to\infty$, uniformly in every closed sector
$|\arg z|\le\pi-\eta$, with $\eta>0$, using the continuation of the
same branch to the slit plane. In particular, for fixed real $s$,
\begin{equation}\label{eq:pre-stirling-phase}
\operatorname{Im}\operatorname{Log}\Gamma(x+\mathrm{i}s)
=s\log x+O(x^{-1}),\qquad x\to\infty.
\end{equation}
Binet's second formula gives the exact identity
\begin{equation}\label{eq:pre-binet}
\operatorname{Log}\Gamma(z)
=\left(z-\frac12\right)\operatorname{Log}z-z
+\frac12\log(2\pi)
+2\int_0^\infty\frac{\arctan(t/z)}{e^{2\pi t}-1}\,\mathrm{d}t,
\qquad \operatorname{Re}z>0.
\end{equation}
Here the complex arctangent is the analytic branch obtained from
its real values for positive $z$. See
\cite[equations~5.9.10 and 5.11.1]{DLMF}. At $z=\mathrm{i}s$, $s>0$,
we use the boundary value
\begin{equation}\label{eq:pre-gamma-boundary}
\operatorname{Log}\Gamma(\mathrm{i}s)
=\operatorname{Log}\Gamma(1+\mathrm{i}s)-\operatorname{Log}(\mathrm{i}s).
\end{equation}

\paragraph{The Gauss hypergeometric function.}
For $c\notin\{0,-1,-2,\ldots\}$, the Gauss hypergeometric function is defined by
\begin{equation}\label{eq:pre-hyper-series}
{}_2F_1(a,b;c;z)
=\sum_{j=0}^\infty\frac{(a)_j(b)_j}{(c)_j\,j!}z^j,
\qquad |z|<1,
\end{equation}
where $(a)_0=1$ and $(a)_j=a(a+1)\cdots(a+j-1)$ for $j\ge1$.
This series is analytic in the unit disk and may be differentiated
term by term there. It solves
\begin{equation}\label{eq:pre-hyper-ode}
{z(1-z)\mathcal F''+[c-(a+b+1)z]\mathcal F'-ab\mathcal F=0.}
\end{equation}
In particular,
\begin{equation}\label{eq:pre-hyper-small}
{}_2F_1(a,b;c;z)=1+\frac{ab}{c}z+O(z^2),\qquad z\to0.
\end{equation}
For $\operatorname{Re}(c-a-b)>0$, Gauss's summation formula is
\begin{equation}\label{eq:pre-gauss}
{}_2F_1(a,b;c;1)
=\frac{\Gamma(c)\Gamma(c-a-b)}{\Gamma(c-a)\Gamma(c-b)}.
\end{equation}
The value at $1$ is also the limit as $z\rightarrow1$ along the real
interval. This convergence is locally uniform in admissible parameters
when the real part of $c-a-b$ is bounded away from zero.
In the case $c=a+b$, with $\operatorname{Re}a>0$ and
$\operatorname{Re}b>0$, we have
\begin{equation}\label{eq:pre-hyper-log}
{}_2F_1(a,b;a+b;z)
=\frac{\Gamma(a+b)}{\Gamma(a)\Gamma(b)}
\bigl[-\log(1-z)+O(1)\bigr],\qquad z\rightarrow1.
\end{equation}
These formulas
follow from \cite[equations~15.4.20 and
15.8.10]{DLMF}.

\subsection[Spectral multipliers and a Riesz mean estimate]{{Spectral multipliers and a Riesz mean estimate}}

We use spherical inversion on $\mathbb H^d_\kappa$ in the normalization
of \cite{Helgason,Strichartz}. The Plancherel density with respect to the spectral parameter $t$ is
\begin{equation}\label{eq:plancherel-density}
w_{d,\kappa}(t)=
\frac{2\kappa^{d-1}}{(4\pi)^{d/2}\Gamma(d/2)}
\left|\frac{\Gamma\left((d-1)/2+\mathrm i t/\kappa\right)}
{\Gamma(\mathrm i t/\kappa)}\right|^2,
\qquad t>0.
\end{equation}
Let $\Phi_{d,\kappa,t}(r)$ be the smooth radial spherical function
normalized by $\Phi_{d,\kappa,t}(0)=1$. It solves
\[
-\Phi''-(d-1)\kappa\coth(\kappa r)\Phi'
=(\rho_{d,\kappa}^2+t^2)\Phi.
\]
{Let $d_\kappa(x,y)$ denote the hyperbolic distance between $x$ and $y$. For a bounded Borel function $f$ with compact support, the kernel of}
$f(H)$ is
\begin{equation}\label{eq:multiplier-kernel}
K_f(x,y)=\int_0^\infty f(\rho_{d,\kappa}^2+t^2)
\Phi_{d,\kappa,t}(d_\kappa(x,y))w_{d,\kappa}(t)\,\mathrm dt.
\end{equation}
\begin{remark}
Formulas~\eqref{eq:plancherel-density} and
\eqref{eq:multiplier-kernel} follow from the Fourier inversion
and Plancherel formulas in \cite[equations~(5.5)--(5.10)]{Strichartz},
together with curvature scaling. In dimension two the density reduces to
\begin{equation}\label{eq:planar-plancherel-density}
	w_{2,\kappa}(t)=\frac{t}{2\pi}\tanh(\pi t/\kappa),
	\qquad t>0,
\end{equation}
consistently with the Mehler--Fock formula and
\cite[equation~(114)]{Marklof}. 	
\end{remark}

We next recall the operator concepts needed for the localized
trace formula. In the following definitions, $\mathcal H$
denotes a separable Hilbert space, and $(e_n)$ denotes an
orthonormal basis of $\mathcal H$.

\begin{definition}\label{def:trace-class}
	A compact operator $T$ on
	$\mathcal H$ is \emph{trace class} if its singular values, namely the
	eigenvalues of $|T|=(T^*T)^{1/2}$ counted with multiplicity, satisfy
	\[
	\sum_{n\ge1}s_n(T)<\infty.
	\]
	For such an operator, its trace is
	\[
	\operatorname{Tr} T=\sum_{n\ge1}\langle Te_n,e_n\rangle_{\mathcal H},
	\]
	where $(e_n)$ is an orthonormal basis of $\mathcal H$.
\end{definition}

The trace series is absolutely convergent and its value is independent
of the choice of orthonormal basis.
In particular, a positive bounded operator $T$ is trace class if and only
if this sum is finite. For a
positive compact operator the trace is the sum of its eigenvalues,
counted with multiplicity.

\begin{definition}\label{def:hilbert-schmidt}
	A bounded operator $S$ on $\mathcal H$ is \emph{Hilbert--Schmidt} if
	\[
	\|S\|_{\mathrm{HS}}^2=\sum_{n\ge1}\|Se_n\|_{\mathcal H}^2<\infty.
	\]
\end{definition}
$\|S\|_{\mathrm{HS}}$ is called the Hilbert--Schmidt norm. This norm is also independent of the choice of orthonormal basis.
For a Hilbert--Schmidt operator $S$, the operator $S^*S$ is positive and
trace class, and
\[
\operatorname{Tr}(S^*S)=\|S\|_{\mathrm{HS}}^2.
\]
If $S$ acts on $L^2(M)$ by an integral kernel $L(x,y)$ in
$L^2(M\times M)$, we have
\[
\|S\|_{\mathrm{HS}}^2
=\int_M\int_M |L(x,y)|^2\,\mathrm{d} V_g(x)\,\mathrm{d} V_g(y).
\]
These facts explain the diagonal kernel formula for the trace used
below. For the standard operator theory, we refer the reader to \cite{SimonTrace}.

\begin{lemma}\label{lem:trace}
{Let $d\ge2$, let $\Omega\subset\mathbb H^d_\kappa$ be measurable with
$V=\operatorname{Vol}_\kappa(\Omega)<\infty$, and let $f$ be a
nonnegative bounded Borel function with compact support.}
Then $\mathbf1_\Omega f(H)\mathbf1_\Omega$ is positive and trace class, and
\begin{equation}\label{eq:trace-formula}
\operatorname{Tr}[\mathbf1_\Omega f(H)\mathbf1_\Omega]
=V\int_0^\infty f(\rho_{d,\kappa}^2+t^2)w_{d,\kappa}(t)\,\mathrm dt.
\end{equation}
Here $\mathbf1_\Omega$ is multiplication by the characteristic
function of $\Omega$ on $L^2(\mathbb H^d_\kappa)$.
\end{lemma}
\begin{proof}
Put $h=\sqrt f$ and $S=h(H)\mathbf1_\Omega$.
Plancherel's identity for \eqref{eq:multiplier-kernel} gives
\[
\int_{\mathbb H^d_\kappa}|K_h(x,y)|^2\,\mathrm dV_\kappa(x)
=\int_0^\infty f(\rho_{d,\kappa}^2+t^2)w_{d,\kappa}(t)\,\mathrm dt.
\]
The right side is finite and independent of $y$.
Integrating over $y\in\Omega$ proves that $S$ is Hilbert--Schmidt.
Since $S^*S=\mathbf1_\Omega f(H)\mathbf1_\Omega$, its trace equals
$\|S\|_{\mathrm{HS}}^2$, which gives \eqref{eq:trace-formula}.
\end{proof}

\begin{lemma}\label{lem:riesz}
For every nonempty open set $\Omega\subset\mathbb H^d_\kappa$ of finite
volume $V$ and every $\lambda\in\mathbb R$,
\begin{equation}\label{eq:riesz-estimate}
R_\Omega(\lambda)
\le V\int_0^\infty
(\lambda-\rho_{d,\kappa}^2-t^2)_+w_{d,\kappa}(t)\,\mathrm dt.
\end{equation}
\end{lemma}
\begin{proof}
Since $\lambda_j(\Omega)>\rho_{d,\kappa}^2$ for every $j\ge1$, if $\lambda\le\rho_{d,\kappa}^2$ then both sides of \eqref{eq:riesz-estimate} vanish. We therefore assume that $\lambda>\rho_{d,\kappa}^2$. Let \(\{u_j\}_{j\ge 1}\) denote the orthonormal basis of \(L^2(\Omega)\) consisting of Dirichlet eigenfunctions. Let $\widetilde u_j$ denote the zero extension of
$u_j$ to $\mathbb H^d_\kappa$. Since $u_j\in H^1_0(\Omega)$, its
zero extension belongs to $H^1(\mathbb H^d_\kappa)$. Moreover,
the family $\{\widetilde u_j\}_{j\ge1}$ is orthonormal in
$L^2(\mathbb H^d_\kappa)$, and $\|\widetilde u_j\|_{L^2(\mathbb H^d_\kappa)}=1$. Zero extension preserves the quadratic form, so $q_H[\widetilde u_j]=q_\Omega[u_j]$. In addition, the eigenvalue equation and normalization give $q_\Omega[u_j]=\lambda_j(\Omega)$.

Set
\[
J_\lambda=\{j\ge1:\lambda_j(\Omega)<\lambda\}.
\]
This is a finite set. Let $P_H$ denote the projection-valued spectral measure of
$H=-\Delta_{\mathbb H^d_\kappa}$. For
$v\in H^1(\mathbb H^d_\kappa)$, define the finite positive measure
\[
\mu_v(B)
=
\langle P_H(B)v,v\rangle_{L^2(\mathbb H^d_\kappa)}
\]
for every Borel set $B\subset\mathbb R$. The spectral theorem gives $\mu_v(\mathbb R)=\|v\|_{L^2(\mathbb H^d_\kappa)}^2$. Since $H$ is nonnegative, the representation of its associated
closed quadratic form yields
\[
D(q_H)=D(H^{1/2})=H^1(\mathbb H^d_\kappa),
\]
and
\[
q_H[v]
=
\int_{\mathbb H^d_\kappa}|\nabla_g v|_g^2\,\mathrm dV_\kappa
=
\|H^{1/2}v\|_{L^2(\mathbb H^d_\kappa)}^2.
\]
Thus, applying the spectral theorem to $H^{1/2}$, we obtain
\[
q_H[v]
=
\|H^{1/2}v\|_{L^2(\mathbb H^d_\kappa)}^2
=
\int_{\mathbb R}\xi\,\mathrm d\mu_v(\xi).
\]
Consequently,
\[
\lambda\|v\|_{L^2(\mathbb H^d_\kappa)}^2-q_H[v]
=
\int_{\mathbb R}(\lambda-\xi)\,\mathrm d\mu_v(\xi).
\]
On the other hand, the functional calculus gives
\[
\langle(\lambda-H)_+v,v\rangle_{L^2(\mathbb H^d_\kappa)}
=
\int_{\mathbb R}(\lambda-\xi)_+\,\mathrm d\mu_v(\xi).
\]
By the scalar inequality $\lambda-\xi\le(\lambda-\xi)_+$, we obtain
\[
\lambda\|v\|_{L^2(\mathbb H^d_\kappa)}^2-q_H[v]
\le
\langle(\lambda-H)_+v,v\rangle_{L^2(\mathbb H^d_\kappa)}
\]
for every $v\in H^1(\mathbb H^d_\kappa)$. Moreover, applying this inequality to $\widetilde u_j$, with
$j\in J_\lambda$, we have
\[
\lambda-\lambda_j(\Omega)
\le
\langle(\lambda-H)_+\widetilde u_j,\widetilde u_j
\rangle_{L^2(\mathbb H^d_\kappa)}.
\]
To estimate the sum of these values, we define
\[
f_\lambda(\xi)
=
(\lambda-\xi)_+
\mathbf1_{[\rho_{d,\kappa}^2,\lambda]}(\xi).
\]
This is a nonnegative bounded Borel function with compact support.
Since the spectrum of $H$ is contained in
$[\rho_{d,\kappa}^2,\infty)$, the functional calculus gives
\[
f_\lambda(H)=(\lambda-H)_+.
\]
Lemma~\ref{lem:trace} therefore shows that $T_\lambda=\mathbf1_\Omega(\lambda-H)_+\mathbf1_\Omega$ is positive and trace class. Since $\mathbf1_\Omega\widetilde u_j=\widetilde u_j$ and
$\mathbf1_\Omega$ is self-adjoint, we have
\[
\langle(\lambda-H)_+\widetilde u_j,\widetilde u_j
\rangle_{L^2(\mathbb H^d_\kappa)}
=
\langle T_\lambda\widetilde u_j,\widetilde u_j
\rangle_{L^2(\mathbb H^d_\kappa)}.
\]
Summing over $j\in J_\lambda$ gives
\[
R_\Omega(\lambda)
\le
\sum_{j\in J_\lambda}
\langle T_\lambda\widetilde u_j,\widetilde u_j
\rangle_{L^2(\mathbb H^d_\kappa)}.
\]
We next extend the finite orthonormal family
$\{\widetilde u_j:j\in J_\lambda\}$ to an orthonormal basis of
$L^2(\mathbb H^d_\kappa)$. Hence
\[
\sum_{j\in J_\lambda}
\langle T_\lambda\widetilde u_j,\widetilde u_j
\rangle_{L^2(\mathbb H^d_\kappa)}
\le
\operatorname{Tr}T_\lambda.
\]
Finally, the trace formula in Lemma~\ref{lem:trace} gives
\[
\operatorname{Tr}T_\lambda
=
V\int_0^\infty
(\lambda-\rho_{d,\kappa}^2-t^2)_+
w_{d,\kappa}(t)\,\mathrm dt.
\]
Combining the preceding inequalities, we prove
\eqref{eq:riesz-estimate}.	
\end{proof}

\section{P\'olya's inequality for hyperbolic balls}\label{sec:balls}

We first consider balls in $\mathbb{H}^d_\kappa$ for any dimension $d\ge2$. Until the final scaling step, we take $\kappa=1$ and set $\rho=(d-1)/2$. As shown in Section~\ref{subsec:singular-oscillation},
separation of variables gives
\[
-w''+\frac{\nu^2-1/4}{\sinh^2r}w=s^2w,
\qquad 0<r<R,
\]
with the Friedrichs condition at $r=0$ and $w(R)=0$.

When $\kappa=1$, {we set $g_{R,s}=G_{1,R,s}$} in \eqref{eq:general-action}. Thus, for real $\nu \ge0$,
\[
g_{R,s}\left(\nu \right)=\frac1\pi\int_0^R
\sqrt{\left(s^2-\nu^2\operatorname{csch}^2r\right)_+}\,\mathrm{d} r.
\]
Let $\nu_{\max}=s\sinh R$ and 
$$
r_\nu =\operatorname{arsinh}\left(\frac{\nu}{s}\right).
$$
If $0\le \nu \le \nu_{\max}$, then $g_{R,s}\left(\nu \right)$ can be represented by
\[
g_{R,s}\left(\nu \right)=\frac1\pi\int_{r_\nu}^R
\sqrt{s^2-\nu^2\operatorname{csch}^2r}\,\mathrm{d} r.
\]
It is easy to see that $g_{R,s}\left(0\right)=sR/\pi$. When $\nu \ge \nu_{\max}$, $g_{R,s}\left(\nu \right)=0$. For $0<\nu <\nu_{\max}$, direct integration gives
\begin{equation}\label{eq:actionclosed}
	\pi g_{R,s}\left(\nu \right)
	=s\operatorname{arcosh}\left(\frac{s\cosh R}{\sqrt{s^2+\nu^2}}\right)-\nu \arctan\left(\frac{\sqrt{s^2\sinh^2R-\nu^2}}{\nu \cosh R}\right).
\end{equation}
Differentiating with respect to $\nu$ then gives
\begin{equation}\label{eq:actionderivative}
	g_{R,s}'\left(\nu \right)=-\frac1\pi\arccos\left(\frac{\nu \coth R}{\sqrt{s^2+\nu^2}}\right)<0.
\end{equation}
A second differentiation gives
\begin{equation}\label{eq:actionsecond}
	g_{R,s}''\left(\nu \right)=\frac{s^2\coth R}
	{\pi\left(s^2+\nu^2\right)\sqrt{s^2-\nu^2\operatorname{csch}^2R}}>0.
\end{equation}
The right derivative at $0$ and the left derivative at $\nu_{\max}$ are \(-1/2\) and $0$, respectively. Hence $g_{R,s}\left(\nu \right)$ is nonnegative, nonincreasing, convex, and $1/2$-Lipschitz. These are precisely the hypotheses of Theorems~\ref{thm:lattice-two} and \ref{thm:lattice-high}.

Let $n_\nu \left(s;R\right)$ be the number of eigenvalues at most $s^2$ of
\eqref{eq:radial} for an arbitrary real $\nu \ge0$. The following proposition provides the key estimate for the
radial eigenvalue counts, which will be combined with the
lattice inequalities to prove P\'olya's inequality for
hyperbolic balls. We first state the result, whose proof is completed in {Subsections~\ref{Jost-Appell}--\ref{phase-to-counting}}.

\begin{proposition}\label{prop:phase}
	For every $s>0$, the quantity in \eqref{eq:epsilon-definition} satisfies $0<\varepsilon\left(s\right)<1/4$. For every $R>0$ and $\nu \ge0$, we have
	\begin{equation}\label{eq:phasebound}
		\left\lfloor g_{R,s}\left(\nu \right)\right\rfloor
		\le n_\nu \left(s;R\right)
		\le\left\lfloor g_{R,s}\left(\nu \right)+\frac14-\varepsilon\left(s\right)\right\rfloor.
	\end{equation}
	In particular, the simpler bound
	\begin{equation}\label{eq:phasebound-basic}
		n_\nu \left(s;R\right)\le\left\lfloor g_{R,s}\left(\nu \right)+\frac14\right\rfloor
	\end{equation}
	also holds.
\end{proposition}

\subsection{The Jost solution and an Appell comparison}\label{Jost-Appell}

{Fix $s>0$ and $\nu \ge0$.} Set $z=e^{-2r}$ and $a_\nu =1/2-\nu$.
We use the special function conventions and identities from
Section~\ref{subsec:special-oscillation}. To estimate the eigenvalue counting function through {Lemma~\ref{lem:pre-oscillation}}, we introduce a complex solution whose
phase determines the zeros of the solution satisfying
the Friedrichs condition at $r=0$.

\paragraph{The Jost solution.}
\begin{definition}[See \cite{KTT}]\label{def:jost}
	For fixed $s>0$ and $\nu\ge0$, the Jost solution
	$f_\nu$ of \eqref{eq:radial} on $(0,\infty)$ is
	characterized by
	\[
	e^{-\mathrm{i}sr}f_\nu(r)\rightarrow1,
	\qquad
	e^{-\mathrm{i}sr}f_\nu'(r)\rightarrow\mathrm{i}s
    \quad \text{as}\quad r\to\infty.
	\]
\end{definition}

For equation \eqref{eq:radial}, we claim that its Jost solution is 
\begin{equation}\label{eq:jost}
	f_\nu \left(r\right)=e^{\mathrm{i} sr}\left(1-z\right)^{a_\nu}
	{}_2F_1\left(a_\nu,a_\nu -\mathrm{i} s;1-\mathrm{i} s;z\right).
\end{equation}
To verify that $f_\nu \left(r\right)$ satisfies \eqref{eq:radial} and is indeed the Jost solution, we set
\[
a=a_\nu =\frac12-\nu,\qquad \gamma=-\frac{\mathrm{i}s}{2}.
\]
Then \(a(a-1)=\nu^2-\frac14\), and
\[
\frac{d^2}{dr^2}
=4z^2\frac{d^2}{dz^2}+4z\frac{d}{dz},
\qquad
\frac{1}{\sinh^2r}=\frac{4z}{(1-z)^2}.
\]
Thus \eqref{eq:radial} is equivalent to
\[
z^2f_{zz}+zf_z+
\left(\frac{s^2}{4}-\frac{a(a-1)z}{(1-z)^2}\right)f=0.
\]
Setting \(f=z^\gamma(1-z)^a{\mathcal F}(z)\) and substituting it into \eqref{eq:radial}, {we obtain}
\[
z^2{\mathcal F}''
+\left((2\gamma+1)z-\frac{2az^2}{1-z}\right){\mathcal F}'
+\left(\gamma^2+\frac{s^2}{4}
-\frac{a(a+2\gamma)z}{1-z}\right){\mathcal F}=0.
\]
Since \(\gamma^2+s^2/4=0\), we have
\[
z(1-z){\mathcal F}''
+\left[1-\mathrm{i}s-(2a+1-\mathrm{i}s)z\right]{\mathcal F}'
-a(a-\mathrm{i}s){\mathcal F}=0.
\]
By \eqref{eq:pre-hyper-ode}, this is the hypergeometric equation satisfied by
\[
{\mathcal F}(z)={}_2F_1(a,a-\mathrm{i}s;1-\mathrm{i}s;z),
\]
which proves that \(f_\nu \) satisfies the radial equation \eqref{eq:radial}.

From the expansion \eqref{eq:pre-hyper-small}, we have
\[
{}_2F_1(a,a-\mathrm{i}s;1-\mathrm{i}s;z)
=1+\frac{a(a-\mathrm{i}s)}{1-\mathrm{i}s}z+O(z^2)\quad \text{as}\quad z\to0.
\]
Combining this with $(1-z)^a=1-az+O(z^2)$, we obtain
\[
(1-z)^a{}_2F_1(a,a-\mathrm{i}s;1-\mathrm{i}s;z)
=1+\frac{a(a-1)}{1-\mathrm{i}s}z+O(z^2)\quad \text{as}\quad z\to0.
\]
Consequently, as $r\to\infty$,
\[
f_\nu (r)=e^{\mathrm{i}sr}
\left(
1+\frac{a(a-1)}{1-\mathrm{i}s}e^{-2r}
+O(e^{-4r})
\right).
\]
In particular,
\begin{equation}\label{asymptotic-formulas-f}
	f_\nu (r)=e^{\mathrm{i}sr}\left(1+O(e^{-2r})\right).
\end{equation}
Differentiating the convergent expansion, we have 
\begin{equation}\label{asymptotic-formulas}
	f_\nu '(r)=\mathrm{i}s\,e^{\mathrm{i}sr}+O(e^{-2r})\quad \text{and}\quad f_\nu ''(r)=-s^2e^{\mathrm{i}sr}+O(e^{-2r}).
\end{equation}

Since \eqref{eq:radial} has real coefficients, both $f_\nu$ and
$\overline{f_\nu}$ are solutions. Their Wronskian is
\[
W_\nu (r)=\overline{f_\nu (r)}f_\nu '(r)
-\overline{f_\nu '(r)}f_\nu (r).
\]
Differentiating and using the radial equation \eqref{eq:radial}, we obtain
\[
W_\nu '(r)
=\overline{f_\nu (r)}f_\nu ''(r)
-\overline{f_\nu ''(r)}f_\nu (r)
=0.
\]
Thus $W_\nu$ is constant.
The asymptotic formulas \eqref{asymptotic-formulas} show that $W_\nu (r)\to2\mathrm{i} s$ as $r\to\infty$.
Consequently, $W_\nu (r)=2\mathrm{i} s$ for every $r>0$.
Since $W_\nu =2\mathrm{i}\operatorname{Im}(\overline{f_\nu}f_\nu ')$,
we obtain
\begin{equation}\label{eq:wronskian}
	\operatorname{Im}\left(\overline{f_\nu}f_\nu '\right)=s>0.
\end{equation}
In particular, $f_\nu\neq0$ on $(0,\infty)$. Therefore, we can write
\[
f_\nu (r)=|f_\nu (r)|e^{\mathrm{i}\theta_\nu (r)}
\]
with a continuous phase $\theta_\nu$ chosen so that
$\theta_\nu (r)-sr\to0$ as $r\to\infty$.
Differentiating this representation yields
\[
\overline{f_\nu (r)}f_\nu '(r)
=|f_\nu (r)|\frac{\mathrm d}{\mathrm dr}|f_\nu (r)|
+\mathrm{i}|f_\nu (r)|^2\theta_\nu '(r).
\]
Taking imaginary parts and applying \eqref{eq:wronskian},
we conclude that
\begin{equation}\label{formula-phase-derivative}
	\theta_\nu '(r)=\frac{s}{|f_\nu (r)|^2}>0.
\end{equation}
Thus the phase is strictly increasing on $(0,\infty)$.

The positivity of $\theta_\nu'$ only shows that the phase
is strictly increasing. To obtain an upper bound for the
eigenvalue counting function, we need to compare the phase with $\pi g_{R,s}(\nu)$.
{For this purpose, we prove a stronger inequality.} For $r>r_\nu$, let 
$$
\zeta =\operatorname{csch}^2r\quad \text{and}\quad \chi =s^2-\nu^2\zeta .
$$ 
Then $\chi >0$. Set 
$$
q=\chi +\frac{\zeta}{4},
$$ 
so the radial equation \eqref{eq:radial} is $f''+qf=0$. We want to compare 
$$
h=\left|f_\nu \right|^2
$$ 
with the explicit function
\[
h_0=s\chi^{-1/2}>0.
\]

\paragraph{The Appell equation.}
For a real-valued function $q\in C^1$, direct differentiation shows that the squared modulus
$h=|f|^2$ of any solution of $f''+qf=0$ satisfies
\[
h'''+4qh'+2q'h=0.
\]
Following \cite[Section~2.2]{Bremer}, we refer to this
equation as \textit{the Appell equation}. The Appell operator is defined by
\[
L_qh=h'''+4qh'+2q'h.
\]
A direct calculation shows that $h_0$ satisfies
$L_{\tilde q}h_0=0$, where
\[
\tilde q=\chi+\frac{\chi''}{4\chi}
-\frac{5(\chi')^2}{16\chi^2}.
\]
Set $T=q-\tilde q$. Using $\zeta ''=4\zeta +6\zeta^2$ and $\left(\zeta '\right)^2=4\zeta^2\left(1+\zeta \right)$ gives
\[
T=\frac \zeta 4+
\frac{\nu^2\zeta \left(4+6\zeta \right)}{4\chi}
+\frac{5\nu^4\zeta^2\left(1+\zeta \right)}{4\chi^2}>0.
\]
It is easy to see that {$T'<0$:} $\zeta$ strictly decreases and $\chi =s^2-\nu^2\zeta$ increases. Also $h_0'\le0$. 

We set $F=L_qh_0$. Since $L_{\tilde q }h_0=0$, we have $F=4Th_0'+2T'h_0<0$. We claim that the Green kernel for $L_q$ is
\begin{equation}\label{eq:kernel}
	K(r,t)=\frac{
		\left[\operatorname{Im}\left(\overline{f_\nu(r)}f_\nu(t)\right)\right]^2
	}{2s^2},
	\qquad t\ge r.
\end{equation}
For completeness, we give the proof of this assertion below. Indeed, for fixed $t$, the function
$r\mapsto\operatorname{Im}(\overline{f_\nu(r)}f_\nu(t))$
solves the radial equation, and hence its square satisfies
the Appell equation. Thus $L_q^{(r)}K=0$, where $(r)$ denotes differentiation with respect to $r$. Similarly, we also obtain
\[
K_{ttt}(r,t)+4q(t)K_t(r,t)+2q'(t)K(r,t)=0.
\]

{To justify the limit below, we next establish the required asymptotic relations at infinity. For $h=|f_\nu|^2$, formulas
\eqref{asymptotic-formulas-f} and \eqref{asymptotic-formulas} give}
\[
h(t)=1+O(e^{-2t}),
\]
\[
h'(t)
=2\operatorname{Re}
\bigl(\overline{f_\nu(t)}f_\nu'(t)\bigr)
=O(e^{-2t}),
\]
\[
h''(t)
=2|f_\nu'(t)|^2
+2\operatorname{Re}
\bigl(\overline{f_\nu(t)}f_\nu''(t)\bigr)
=O(e^{-2t}).
\]
Since $\operatorname{csch}^2t$ and its first two
derivatives are $O(e^{-2t})$, we have
\[
\chi(t)=s^2+O(e^{-2t}),
\qquad
\chi'(t)=O(e^{-2t}),
\qquad
\chi''(t)=O(e^{-2t}).
\]
It follows that $h_0=s\chi^{-1/2}$ satisfies
\[
{h_0(t)=1+O(e^{-2t}),\qquad h_0'(t)=O(e^{-2t}),\qquad h_0''(t)=O(e^{-2t}).}
\]
Therefore, setting $w=h_0-h$, we obtain
\[
w^{(j)}(t)=O(e^{-2t}),
\qquad j=0,1,2,
\qquad t\to\infty.
\]
It is also easy to see that
$$
L_qw=F.
$$
Next, fix $r>r_\nu$ and set
\[
a_r(t)=\operatorname{Im}
\bigl(\overline{f_\nu(r)}f_\nu(t)\bigr).
\]
Then $a_r''(t)=-q(t)a_r(t)$, and
\[
K(r,t)=\frac{a_r(t)^2}{2s^2},\quad K_t(r,t)=\frac{a_r(t)a_r'(t)}{s^2},\quad 
K_{tt}(r,t)=\frac{a_r'(t)^2-q(t)a_r(t)^2}{s^2}.
\]
{The asymptotic formulas \eqref{asymptotic-formulas-f} and \eqref{asymptotic-formulas}} show that $f_\nu(t)$ and
$f_\nu'(t)$ remain bounded as $t\to\infty$.
Since $q(t)\to s^2$, these identities imply that the functions $K(r,t)$,
$K_t(r,t)$, and $K_{tt}(r,t)$ remain bounded as
$t\to\infty$. At $t=r$, the Wronskian identity gives
\[
a_r(r)=0,
\qquad
a_r'(r)
=\operatorname{Im}
\bigl(\overline{f_\nu(r)}f_\nu'(r)\bigr)=s.
\]
Hence
\begin{equation}\label{diagonal-identities}
	K(r,r)=K_t(r,r)=0,
	\qquad
	K_{tt}(r,r)=1.
\end{equation}
Multiplying both sides of $L_qw=F$ by $K(r,t)$
and integrating by parts over $t\in[r,U]$, we obtain	
\begin{align*}
	\int_r^U K(r,t)F(t)\,\mathrm{d}t
	&=\int_r^U K(r,t)
	\bigl(w'''(t)+4q(t)w'(t)+2q'(t)w(t)\bigr)
	\,\mathrm{d}t\\
	&=\Bigl[
	K(r,t)w''(t)-K_t(r,t)w'(t)\\
	&\hspace{2.5em}
	+\bigl(K_{tt}(r,t)+4q(t)K(r,t)\bigr)w(t)
	\Bigr]_{t=r}^{t=U}\\
	&\quad-\int_r^U
	\bigl(K_{ttt}(r,t)+4q(t)K_t(r,t)
	+2q'(t)K(r,t)\bigr)w(t)\,\mathrm{d}t.
\end{align*}
The last integral vanishes by the Appell equation for
$K(r,\cdot)$. The boundary expression at $t=r$ equals
$w(r)$ by \eqref{diagonal-identities}. Thus
\[
\begin{aligned}
	w(r)
	&=K(r,U)w''(U)-K_t(r,U)w'(U)\\
	&\quad+
	\bigl(K_{tt}(r,U)+4q(U)K(r,U)\bigr)w(U)
	-\int_r^U K(r,t)F(t)\,\mathrm{d}t.
\end{aligned}
\]
The boundary terms at $U$ tend to zero as
$U\to\infty$ because their coefficients remain bounded
and $w$, $w'$, and $w''$ tend to zero. Thus, we obtain 
\[
h_0(r)-h(r)
=-\int_r^\infty K(r,t)F(t)\,\mathrm{d}t>0.
\]	
The integral converges by the exponential decay of $F$. The inequality is strict because $F<0$, $K\ge0$, and
{$K(r,t)>0$ for $t>r$ sufficiently close to $r$.}

Thus, we have proved that
\[
\left|f_\nu \left(r\right)\right|^2<\frac{s}{\sqrt{s^2-\nu^2\operatorname{csch}^2r}},\qquad r>r_\nu .
\]
Using \eqref{formula-phase-derivative}, we obtain
\begin{equation}\label{eq:modulus}
\theta_\nu '\left(r\right)>\sqrt{s^2-\nu^2\operatorname{csch}^2r},\qquad r>r_\nu.
\end{equation}

\subsection{Phase estimate}\label{phase-estimate}

{In this subsection, we use the branch} of $\operatorname{Log}\Gamma$ fixed in
Section~\ref{subsec:special-oscillation}. Set
\[
\alpha_\nu \left(s\right)=\operatorname{Im}\operatorname{Log}\Gamma\left(1+\mathrm{i} s\right)
-\operatorname{Im}\operatorname{Log}\Gamma\left(\nu +1/2+\mathrm{i} s\right).
\]
For $\nu >0$, Gauss summation \eqref{eq:pre-gauss} applied to \eqref{eq:jost} gives
\begin{equation}\label{eq:originpositive}
	f_\nu \left(r\right)\sim \left(2r\right)^{1/2-\nu}
	\frac{\Gamma\left(1-\mathrm{i} s\right)\Gamma\left(2\nu \right)}
	{\Gamma\left(\nu +1/2\right)\Gamma\left(\nu +1/2-\mathrm{i} s\right)}.
\end{equation}
For $\nu =0$, \eqref{eq:pre-hyper-log} gives
\begin{equation}\label{eq:originzero}
	f_0\left(r\right)\sim -\sqrt{2r}\log\left(2r\right)
	\frac{\Gamma\left(1-\mathrm{i} s\right)}{\Gamma\left(1/2\right)\Gamma\left(1/2-\mathrm{i} s\right)}.
\end{equation}

For \(\nu \ge 0\), we wish to obtain the phase of $f_\nu \left(r\right)$ at \(r=0\). We begin with the case \(\nu>0\). In this case, $\Gamma(2\nu)$, $\Gamma(\nu+1/2)$ and $(2r)^{1/2-\nu}$ are positive. Hence the phase of the
leading coefficient in \eqref{eq:originpositive} is
determined by
\[
\frac{\Gamma(1-\mathrm{i}s)}
{\Gamma(\nu+1/2-\mathrm{i}s)}.
\]
{Using $\Gamma(\overline z)=\overline{\Gamma(z)}$, we obtain
\[
\begin{aligned}
	\operatorname{Im}\operatorname{Log}\Gamma(1-\mathrm{i}s)
	-\operatorname{Im}\operatorname{Log}\Gamma
	\left(\nu+\frac12-\mathrm{i}s\right)=-\alpha_\nu(s).
\end{aligned}
\]}
Hence
\[
\frac{\Gamma(1-\mathrm{i}s)}
{\Gamma(\nu+1/2-\mathrm{i}s)}
=
\left|
\frac{\Gamma(1-\mathrm{i}s)}
{\Gamma(\nu+1/2-\mathrm{i}s)}
\right|
e^{-\mathrm{i}\alpha_\nu(s)}.
\]
Thus, we have
\[
e^{\mathrm{i}\theta_\nu(0)}
=e^{-\mathrm{i}\alpha_\nu(s)}.
\]
Consequently,
\[
\theta_\nu(0)+\alpha_\nu(s)=2\pi k(\nu),
\qquad k(\nu)\in\mathbb Z.
\]
The origin asymptotics are locally uniform in $\nu>0$
and the normalized phase depends continuously on $\nu$.
{Thus $k(\nu)$ is continuous and integer-valued on $(0,\infty)$, and hence constant.} To determine this constant, we take $\nu=1/2$. The Jost solution is $f_{1/2}(r)=e^{\mathrm{i}sr}$. Its normalized phase is $\theta_{1/2}(r)=sr$, so
$\theta_{1/2}(0)=0$. Moreover, the definition of
$\alpha_\nu(s)$ gives $\alpha_{1/2}(s)=0$. Hence $k(1/2)=0$. Therefore $k(\nu)=0$ for all $\nu>0$.
We conclude that
\begin{equation}\label{eq:originphase}
	\theta_\nu(0)=-\alpha_\nu(s),
	\qquad \nu>0.
\end{equation}

For the case $\nu =0$, we next establish a Gamma inequality, independently of \eqref{eq:originphase}. We define
\[
C_\nu \left(s\right)=\int_0^\nu \arctan\left(\frac{s}{t}\right)\,\mathrm{d} t.
\]
For $\nu >0$, direct integration gives
\[
C_\nu \left(s\right)=\nu \arctan\left(\frac{s}{\nu} \right)+\frac{s}{2}\log\left(1+\frac{\nu^2}{s^2}\right).
\]
Set
\[
D_\nu \left(s\right)=\frac\pi4-C_\nu \left(s\right)-\alpha_\nu \left(s\right),
\]
{with the continuous value $C_0(s)=0$.}  We set 
\[
a(t)=\arctan\left(\frac{s}{t}\right)
\]
and \(a(0)=\pi/2\). By the recurrence formula for the Gamma function, we have
\[
\begin{aligned}
	D_{\nu+1}(s)-D_\nu(s)
	&=\alpha_\nu(s)-\alpha_{\nu+1}(s)
	-\int_\nu^{\nu+1}a(t)\,\mathrm{d}t\\
	&=\operatorname{Im}\left[
	\operatorname{Log}\Gamma(\nu+3/2+\mathrm{i}s)
	-\operatorname{Log}\Gamma(\nu+1/2+\mathrm{i}s)
	\right]
	-\int_\nu^{\nu+1}a(t)\,\mathrm{d}t\\
	&=\operatorname{Im}\operatorname{Log}(\nu+1/2+\mathrm{i}s)
	-\int_\nu^{\nu+1}a(t)\,\mathrm{d}t\\
	&=\arctan\left(\frac{s}{\nu+1/2}\right)
	-\int_\nu^{\nu+1}a(t)\,\mathrm{d}t\\
	&=a(\nu+1/2)-\int_\nu^{\nu+1}a(t)\,\mathrm{d}t<0,
\end{aligned}
\]
where the last inequality follows from strict convexity of $a$. Let $x=\nu+j$. Then
\[
C_x(s)
=x\arctan(s/x)+\frac{s}{2}\log(1+x^2/s^2)
=s+s\log x-s\log s+o(1)
\]
as $x\to\infty$. Stirling's formula
\eqref{eq:pre-stirling-phase} yields
\[
\operatorname{Im}\operatorname{Log}
\Gamma(x+1/2+\mathrm{i}s)
=s\log x+o(1).
\]
Substituting these asymptotics into the definition of $D_x(s)$,
we obtain
\[
\begin{aligned}
	D_x(s)
	&=\frac{\pi}{4}-C_x(s)
	-\operatorname{Im}\operatorname{Log}\Gamma(1+\mathrm{i}s)
	+\operatorname{Im}\operatorname{Log}
	\Gamma(x+1/2+\mathrm{i}s)\\
	&=\frac{\pi}{4}+s\log s-s
	-\operatorname{Im}\operatorname{Log}\Gamma(1+\mathrm{i}s)
	+o(1).
\end{aligned}
\]
Using the recurrence relation
\[
\operatorname{Log}\Gamma(\mathrm{i}s)
=\operatorname{Log}\Gamma(1+\mathrm{i}s)
-\operatorname{Log}(\mathrm{i}s),
\]
where $\operatorname{Log}(\mathrm{i}s)=\log s+\mathrm{i}\pi/2$,
we conclude that
\[
\lim_{j\to\infty}D_{\nu+j}(s)
=\mathcal T(s)
:=s\log s-s
-\operatorname{Im}\operatorname{Log}\Gamma(\mathrm{i}s)
-\frac{\pi}{4}.
\]
Here $\operatorname{Log}\Gamma(\mathrm{i}s)$ is understood
as the boundary value from the right half-plane. Binet's second formula \eqref{eq:pre-binet} gives
\[
\mathcal T\left(s\right)=\int_0^\infty
\frac{\log\left|\dfrac{s+t}{s-t}\right|}{e^{2\pi t}-1}\,\mathrm{d} t>0.
\]
In the notation of \eqref{eq:epsilon-definition}, $\mathcal T\left(s\right)=\pi\varepsilon\left(s\right)$. We obtain this identity by taking imaginary parts of Binet's
integral at $z=\delta+\mathrm{i} s$ and letting $\delta\rightarrow0$.
For $\delta>0$, its relevant imaginary part is
\[
\begin{aligned}
	\operatorname{Im}\left(
	2\arctan\left(\frac{t}{\delta+\mathrm{i}s}\right)
	\right)
	&=\operatorname{Im}\left\{
	-\mathrm{i}\left[
	\operatorname{Log}\bigl(\delta+\mathrm{i}(s+t)\bigr)
	-\operatorname{Log}\bigl(\delta+\mathrm{i}(s-t)\bigr)
	\right]\right\}\\
	&=-\log\left|\delta+\mathrm{i}(s+t)\right|
	+\log\left|\delta+\mathrm{i}(s-t)\right|\\
	&=-\frac12\log\bigl(\delta^2+(s+t)^2\bigr)
	+\frac12\log\bigl(\delta^2+(s-t)^2\bigr)\\
	&=-\frac12\log\left(
	\frac{\delta^2+(s+t)^2}{\delta^2+(s-t)^2}
	\right).
\end{aligned}
\]
For $t\ne s$, this tends to $-\log\left|\left(s+t\right)/\left(s-t\right)\right|$.  The singularity at $t=s$ is logarithmic and
integrable. The integrand is bounded near zero and it decays
exponentially at infinity.  These bounds justify the limit. It follows that
\begin{equation}\label{eq:Dpositive}
	D_\nu \left(s\right)>\mathcal T\left(s\right)>0\qquad\left(\nu \ge0\right).
\end{equation}

We can now verify \eqref{eq:originphase} at $\nu =0$. For $\nu>0$, the monotonicity of the phase yields
\[
\theta_\nu(r)>\theta_\nu(0)=-\alpha_\nu(s),
\qquad r>0.
\]
For fixed $r>0$, letting $\nu\rightarrow0$ and using
continuity in $\nu$, we obtain
$\theta_0(r)\ge-\alpha_0(s)$. Hence
\[
\theta_0(0)\ge-\alpha_0(s).
\]
On the other hand, \eqref{eq:modulus} gives
$\theta_0'(r)>s$ on $(0,\infty)$. Thus
$\theta_0(r)-sr$ is strictly increasing and tends to zero
as $r\to\infty$. For any fixed $R>0$, it follows that
\[
\theta_0(0)
=\lim_{r\rightarrow0}\bigl(\theta_0(r)-sr\bigr)
\le\theta_0(R)-sR<0.
\]
Combining these inequalities yields
\[
-\alpha_0(s)\le\theta_0(0)<0.
\]
By \eqref{eq:originzero}, there is an integer $j$ such that
\[
\theta_0(0)=-\alpha_0(s)+2\pi j.
\]
The lower bound gives $j\ge0$ and the upper bound and
\eqref{eq:Dpositive} imply
\[
2\pi j<\alpha_0(s)<\frac{\pi}{4}.
\]
Therefore $j=0$ and
\[
\theta_0(0)=-\alpha_0(s),
\qquad 0<\alpha_0(s)<\frac{\pi}{4}.
\]
Thus \eqref{eq:originphase} holds for every $\nu\ge0$.

Since $D_0\left(s\right)=\pi/4-\alpha_0\left(s\right)$, inequality \eqref{eq:Dpositive} also gives
\[
0<\mathcal T\left(s\right)<\frac\pi4.
\]
This proves the asserted bounds on $\varepsilon\left(s\right)$.

\begin{remark}
Passing to the limit $\nu\rightarrow0$ in
\eqref{eq:originphase} would require
\[
\lim_{\nu\rightarrow0}\lim_{r\rightarrow0}\theta_\nu(r)
=
\lim_{r\rightarrow0}\lim_{\nu\rightarrow0}\theta_\nu(r).
\]
Continuity in $\nu$ for each fixed $r>0$ alone does not
justify this interchange. Our argument avoids the interchange of the two limits.
\end{remark}

\subsection{From phase to eigenvalue counting}\label{phase-to-counting}

For $r\ge r_\nu$, we define
\[
S_\nu \left(r\right)=\int_{r_\nu}^r\sqrt{s^2-\nu^2\operatorname{csch}^2t}\,\mathrm{d} t.
\]
Set $S_\nu \left(r\right)=0$ for $0\le r\le r_\nu$.
{From \eqref{eq:actionclosed}, we obtain}
\[
S_\nu \left(r\right)=sr-C_\nu \left(s\right)+o\left(1\right)\qquad r\to\infty.
\]
For $r>r_\nu$, let 
$$
{e_\nu \left(r\right)=C_\nu(s)+S_\nu} \left(r\right)-\theta_\nu \left(r\right).
$$
By \eqref{eq:modulus}, $e_\nu '<0$ and $\lim_{r\to\infty}e_\nu(r)=0$. Thus $e_\nu \left(r\right)>0$. Since $S_\nu \left(r_\nu \right)=0$, we have
\[
{e_\nu \left(r_\nu \right)=C_\nu(s)-\theta_\nu \left(r_\nu \right)\le C_\nu(s)+\alpha_\nu(s).}
\]
By the definition of $D_\nu$,
\[
{C_\nu(s)+\alpha_\nu(s)=\frac\pi4-D_\nu(s).}
\]
Consequently, due to \eqref{eq:Dpositive}, we have
\begin{equation}\label{eq:turning}
	0<e_\nu \left(r_\nu \right)<\pi\left(\frac14-\varepsilon\left(s\right)\right).
\end{equation}
For $\nu =0$, the value at $r_\nu =0$ is understood as the finite limit from the right.

{We claim that the real solution of the radial equation \eqref{eq:radial} satisfying the Friedrichs condition at the origin is a nonzero real multiple of}
\[
\operatorname{Im}\left(
e^{-\mathrm{i}\theta_\nu(0)}f_\nu(r)
\right).
\]
To see this, we choose two real solutions $\phi_\nu$
and $\psi_\nu$, where $\phi_\nu$ satisfies the Friedrichs
condition and $\psi_\nu$ is the singular solution.
For $\nu>0$, we normalize them by
\[
\phi_\nu(r)=r^{1/2+\nu}(1+O(r^2)),
\qquad
\psi_\nu(r)=r^{1/2-\nu}(1+o(1)).
\]
For $\nu=0$, they are
\[
\phi_0(r)=\sqrt r(1+O(r^2)),
\qquad
\psi_0(r)=-\sqrt r\log r(1+o(1)).
\]
{We set
\[
f_\nu(r)=A_\nu\psi_\nu(r)+B_\nu\phi_\nu(r),
\qquad A_\nu,B_\nu\in\mathbb C.
\]
Since \(\phi_\nu(r)/\psi_\nu(r)\to0\) as \(r\to0\), we deduce from the {asymptotic behavior} near \(r=0\) that
\[
\frac{A_\nu}{|A_\nu|}
=e^{\mathrm{i}\theta_\nu(0)}.
\]
Consequently,
\[
e^{-\mathrm{i}\theta_\nu(0)}f_\nu(r)
=
|A_\nu|\psi_\nu(r)
+e^{-\mathrm{i}\theta_\nu(0)}B_\nu\phi_\nu(r).
\]
The coefficient of $\psi_\nu$ is now positive and real.
Since $\psi_\nu$ itself is real, we have
\[
\operatorname{Im}\left(
e^{-\mathrm{i}\theta_\nu(0)}f_\nu(r)
\right)
=
\operatorname{Im}\left(
e^{-\mathrm{i}\theta_\nu(0)}B_\nu
\right)\phi_\nu(r).
\]
}This argument also applies when $\psi_\nu$ contains
logarithmic terms. The remaining coefficient is nonzero; otherwise, the rotated Jost solution would be real, contradicting \eqref{eq:wronskian}.

Setting
\[
\beta_\nu(r)=\theta_\nu(r)-\theta_\nu(0),
\]
we obtain a regular solution proportional to $|f_\nu(r)|\sin\beta_\nu(r)$. Since $|f_\nu(r)|>0$, its positive zeros correspond precisely
to $\beta_\nu(r)=k\pi$, with $k=1,2,\ldots$.

For $r>r_\nu$, the definitions of $D_\nu$ and $e_\nu$ give
\begin{equation}\label{eq:phase-excess-identity}
	\beta_\nu \left(r\right)-S_\nu \left(r\right)
	=\frac\pi4-D_\nu \left(s\right)-e_\nu \left(r\right).
\end{equation}
The left side is positive because strict decrease of $e_\nu$ makes it larger than
$\theta_\nu \left(r_\nu \right)-\theta_\nu \left(0\right)\ge0$.
On the other hand, $e_\nu \left(r\right)>0$ and \eqref{eq:Dpositive} imply
\begin{equation}\label{eq:phase-excess-refined}
	0<\beta_\nu \left(r\right)-S_\nu \left(r\right)
	<\pi\left(\frac14-\varepsilon\left(s\right)\right),\qquad r>r_\nu .
\end{equation}
At $r=r_\nu$, we have
\[
\beta_\nu \left(r_\nu \right)=\frac\pi4-D_\nu \left(s\right)-e_\nu \left(r_\nu \right).
\]
Since $\beta_\nu$ is increasing, $\beta_\nu \left(r\right)\le\pi\left(1/4-\varepsilon\left(s\right)\right)$ on $0\le r\le r_\nu$. {Therefore, the regular radial solution has no zeros in $(0,r_\nu]$.} Lemma~\ref{lem:pre-oscillation}, in the phase form
\eqref{eq:pre-oscillation-phase}, gives
\[
n_\nu \left(s;R\right)=\left\lfloor \frac{\beta_\nu \left(R\right)}{\pi}\right\rfloor.
\]
If $R>r_\nu$, using $S_\nu \left(R\right)=\pi g_{R,s}\left(\nu \right)$ in \eqref{eq:phase-excess-refined} and taking integer parts, we have \eqref{eq:phasebound}. If $R\le r_\nu$, it is trivial because $n_\nu \left(s;R\right)=g_{R,s}\left(\nu \right)=0$ and $0<\varepsilon\left(s\right)<1/4$. Finally, \eqref{eq:phasebound-basic} follows by discarding $\varepsilon\left(s\right)$.

\subsection[Proofs of Theorem~\ref{thm:balls} and its corollaries]{{Proofs of Theorem~\ref{thm:balls} and its corollaries}}

We now combine Proposition~\ref{prop:phase} with Theorems~\ref{thm:lattice-two} and \ref{thm:lattice-high} and prove Theorem~\ref{thm:balls}. 

\begin{proof}[Proof of Theorem~\ref{thm:balls}]
First let $\kappa=1$. Separation of variables gives
\[
N_{B_R}\left(\rho^2+s^2\right)
=\sum_{\ell\ge0}{m_{d,\ell}}\,n_{\ell+\left(d-2\right)/2}\left(s;R\right).
\]
Indeed, let $\{Y_{\ell,m}\}_{m=1}^{{m_{d,\ell}}}$ be an orthonormal
basis of spherical harmonics of degree $\ell$. Expanding
\[
u(r,\omega)
=(\sinh r)^{-\rho}
\sum_{\ell=0}^{\infty}\sum_{m=1}^{{m_{d,\ell}}}
w_{\ell,m}(r)Y_{\ell,m}(\omega),
\]
and using the orthogonality of spherical harmonics, we
obtain the radial equation \eqref{eq:radial} for each $w_{\ell,m}$, where $\nu=\ell+(d-2)/2$ and $s^2=\lambda-\rho^2$. Thus each radial eigenvalue \(\mu_{\nu,j}\) yields an eigenvalue \(\rho^2+\mu_{\nu,j}\) for the ball, with \({m_{d,\ell}}\) linearly independent eigenfunctions. The completeness of the spherical
harmonics ensures that all eigenvalues are obtained in this
way and multiplicities are added. Consequently,
\[
\begin{aligned}
	N_{B_R}(\rho^2+s^2)
	&=\sum_{\ell=0}^{\infty}{m_{d,\ell}}\,
	\#\{j:\mu_{\ell+(d-2)/2,j}\le s^2\}\\
	&=\sum_{\ell=0}^{\infty}{m_{d,\ell}}\,
	n_{\ell+(d-2)/2}(s;R).
\end{aligned}
\]

Let
\begin{equation}\label{eq:truncated-action}
	\widetilde g_{R,s}\left(\nu \right)
	=\left(g_{R,s}\left(\nu \right)-\varepsilon\left(s\right)\right)_+.
\end{equation}
Since $0<\varepsilon\left(s\right)<1/4$,
\begin{equation}\label{eq:floor-truncation}
	\left\lfloor g_{R,s}\left(\nu \right)+\frac14-\varepsilon\left(s\right)\right\rfloor
	=\left\lfloor \widetilde g_{R,s}\left(\nu \right)+\frac14\right\rfloor.
\end{equation}
Proposition~\ref{prop:phase} therefore gives
\[
N_{B_R}\left(\rho^2+s^2\right)
\le\sum_{\ell\ge0}{m_{d,\ell}}
\left\lfloor \widetilde g_{R,s}\left(\ell+\frac{d-2}{2}\right)+\frac14\right\rfloor.
\]
All sums are finite because $g_{R,s}\left(\nu \right)=0$ for $\nu \ge s\sinh R$.

If $g_{R,s}\left(0\right)\le\varepsilon\left(s\right)$, then $\widetilde g_{R,s}$ is identically zero and so is the counting function. Otherwise, strict decrease of $g_{R,s}$ on $\left(0,s\sinh R\right)$ gives a unique $\nu_\varepsilon$ such that $g_{R,s}\left(\nu_\varepsilon\right)=\varepsilon\left(s\right)$. {$\widetilde g_{R,s}\left(\nu\right)>0$ on $[0,\nu_\varepsilon)$ and $\widetilde g_{R,s}\left(\nu\right)=0$ on $[\nu_\varepsilon,\infty)$.} In this case, it is easy to see that $\widetilde g_{R,s}$ is nonnegative, nonincreasing, convex, and $1/2$-Lipschitz. Applying Theorem~\ref{thm:lattice-two} in dimension two and Theorem~\ref{thm:lattice-high} in higher dimensions, we obtain
\begin{equation}\label{eq:action-integral}
	N_{B_R}\left(\rho^2+s^2\right)
	\le\frac2{\left(d-2\right)!}\int_0^\infty
	\nu^{d-2}\widetilde g_{R,s}\left(\nu \right)\,\mathrm{d} \nu .
\end{equation}
This inequality is strict if $\widetilde g_{R,s}$ is not identically zero.

To evaluate the difference from the Weyl term, we consider the following integral:
\[
\mathcal I_{d,R,s}=\frac2{\left(d-2\right)!}\int_0^\infty \nu^{d-2}g_{R,s}\left(\nu \right)\,\mathrm{d} \nu .
\]
By direct calculation, we have
\[
\begin{aligned}
	\mathcal I_{d,R,s}
	&=\frac{2}{\pi\Gamma(d-1)}
	\int_0^\infty\int_0^R
	\nu^{d-2}
	\sqrt{\left(s^2-\frac{\nu^2}{\sinh^2r}\right)_+}
	\,\mathrm{d}r\,\mathrm{d}\nu\\
	&=\frac{2}{\pi\Gamma(d-1)}
	\int_0^R\int_0^{s\sinh r}
	\nu^{d-2}\sqrt{s^2-\frac{\nu^2}{\sinh^2r}}
	\,\mathrm{d}\nu\,\mathrm{d}r\\
	&=\frac{2s^d}{\pi\Gamma(d-1)}
	\int_0^R\sinh^{d-1}r\,\mathrm{d}r
	\int_0^1t^{d-2}\sqrt{1-t^2}\,\mathrm{d}t.
\end{aligned}
\]
The second equality follows from Tonelli's theorem and
the last uses the substitution $\nu=s\sinh r\,t$. Using \eqref{eq:pre-beta}, we get the following identity:
\[
\int_0^1t^{d-2}\sqrt{1-t^2}\,\mathrm{d} t
=\frac12\mathrm B\left(\frac{d-1}{2},\frac32\right).
\]
The identity \eqref{eq:pre-beta-gamma} and the
duplication formula \eqref{eq:pre-duplication} yield
\[
\mathcal I_{d,R,s}=\frac{s^d}{2^{d-1}\Gamma\left(d/2\right)\Gamma\left(d/2+1\right)}
\int_0^R\sinh^{d-1}r\,\mathrm{d} r.
\]
Using \eqref{eq:ball-volume} and the definition of $C_d$, we conclude that
\begin{equation}\label{eq:volume-integral}
	\mathcal I_{d,R,s}=C_d\operatorname{Vol}_1\left(B_R\right)s^d.
\end{equation}
The identity $g_{R,s}\left(\nu \right)-\widetilde g_{R,s}\left(\nu \right)
=\min\left\{g_{R,s}\left(\nu \right),\varepsilon\left(s\right)\right\}$ now gives
\begin{equation}\label{eq:remainder-identity}
	\frac2{\left(d-2\right)!}\int_0^\infty
	\nu^{d-2}\widetilde g_{R,s}\left(\nu \right)\,\mathrm{d} \nu 
	=C_d\operatorname{Vol}_1\left(B_R\right)s^d-\mathcal D_{d,1}\left(R,s\right).
\end{equation}
The remainder is strictly positive because both $g_{R,s}\left(\nu \right)$ and $\varepsilon\left(s\right)$ are positive on $0<\nu <s\sinh R$. {Equations~\eqref{eq:action-integral} and \eqref{eq:remainder-identity} prove} \eqref{eq:ball-count-refined} at $\kappa=1$. Positivity of the remainder also proves \eqref{eq:ball-count}.

For general \(\kappa\), we apply the result for \(\kappa=1\) at radius \(\kappa R\) and shifted parameter \(s/\kappa\). The substitution $t=\kappa r$ in \eqref{eq:general-action} gives
\[
G_{\kappa,R,s}\left(\nu \right)=g_{\kappa R,s/\kappa}\left(\nu \right).
\]
Consequently,
\[
\mathcal D_{d,\kappa}\left(R,s\right)
=\mathcal D_{d,1}\left(\kappa R,s/\kappa\right).
\]
Eigenvalues scale by a factor \(\kappa^2\) and volume by \(\kappa^{-d}\), so the corresponding factors cancel in the Weyl term. This proves the full counting statements. At $s^2=\lambda_k\left(B_R\right)-\rho_{d,\kappa}^2>0$, the count is at least $k$; rearrangement of \eqref{eq:ball-count} proves \eqref{eq:ball-eigen}.

We finally prove the optimality of the leading coefficient
and the constant shift. Set
\[
{A_d=\frac{(2\pi)^2}{\omega_d^{2/d}}.}
\]
For any fixed hyperbolic ball $B_R$, Weyl's law gives
\[
\lim_{k\to\infty}
\frac{\lambda_k(B_R)}
{\left(k/\operatorname{Vol}_{\kappa}(B_R)\right)^{2/d}}
=A_d.
\]
Subtracting the fixed constant $\rho_{d,\kappa}^2$ does not
change this limit. Hence an inequality
\[
\lambda_k(B_R)\ge
\rho_{d,\kappa}^2+
A\left(\frac{k}{\operatorname{Vol}_{\kappa}(B_R)}\right)^{2/d}
\]
valid for every $k$ would imply $A\le A_d$.
Thus the leading coefficient cannot be increased.

To prove optimality of the shift $\rho_{d,\kappa}^2$, we recall that
\[
\inf\sigma(-\Delta_{\mathbb H_\kappa^d})
=\rho_{d,\kappa}^2.
\]
The variational principle therefore gives
$\lambda_1(B_R)\ge\rho_{d,\kappa}^2$.
Conversely, for every $\varepsilon>0$, the variational
characterization of the spectral bottom provides a nonzero
$u\in C_c^\infty(\mathbb H_\kappa^d)$ such that
\[
\frac{\int_{\mathbb H_\kappa^d}{|\nabla_g u|_g^2}\,
	\mathrm{d}V_\kappa}
{\int_{\mathbb H_\kappa^d}|u|^2\,
	\mathrm{d}V_\kappa}
<\rho_{d,\kappa}^2+\varepsilon.
\]
For all sufficiently large $R$, the support of $u$ lies
in $B_R$, so $\lambda_1(B_R)<\rho_{d,\kappa}^2+\varepsilon$.
With domain monotonicity, this proves
\[
\lambda_1(B_R)\rightarrow\rho_{d,\kappa}^2
\qquad R\to\infty;
\]
see also \cite{Davies,ILW}. If $\rho_{d,\kappa}^2$ could be replaced by a constant
$c>\rho_{d,\kappa}^2$ uniformly over all radii, then the
case $k=1$ would give
\[
\lambda_1(B_R)\ge
c+A_d{\operatorname{Vol}_{\kappa}(B_R)}^{-2/d}.
\]
Since ${\operatorname{Vol}_{\kappa}(B_R)}\to\infty$, letting
$R\to\infty$ would yield $\rho_{d,\kappa}^2\ge c$,
a contradiction. 	
\end{proof}

\begin{corollary}[Berezin--Li--Yau inequality for balls]\label{cor:ball-bly}
For every geodesic ball $B_R\subset\mathbb H^d_\kappa$, $d\ge2$,
and every $s>0$,
\begin{equation}\label{eq:ball-bly-riesz}
R_{B_R}(\rho_{d,\kappa}^2+s^2)
<L^{\mathrm{cl}}_{1,d}\operatorname{Vol}_\kappa(B_R)s^{d+2}.
\end{equation}
For every $m\in\mathbb N$,
\begin{equation}\label{eq:ball-bly-sum}
\sum_{j=1}^m\lambda_j(B_R)
>\rho_{d,\kappa}^2m+
\frac{d}{d+2}\,4\pi^2
\left(\frac{m}{\omega_d\operatorname{Vol}_\kappa(B_R)}\right)^{2/d}m.
\end{equation}
Together with Theorem~\ref{thm:balls}, this shows that balls satisfy
both Conjectures~\ref{conj:laptev-polya} and~\ref{conj:laptev-bly}. Both the leading coefficient and the constant
shift in \eqref{eq:ball-bly-sum} are individually optimal uniformly
over all radii and indices.
\end{corollary}
\begin{proof}
By the definition of the first Riesz mean,
\[
R_{B_R}(\rho_{d,\kappa}^2+s^2)
=\sum_{j=1}^{\infty}
\bigl(\rho_{d,\kappa}^2+s^2-\lambda_j(B_R)\bigr)_+.
\]
Since $\lambda_j(B_R)>\rho_{d,\kappa}^2$, each summand can be
written as
\[
\bigl(\rho_{d,\kappa}^2+s^2-\lambda_j(B_R)\bigr)_+
=\int_0^{s^2}
\mathbf{1}_{\{\lambda_j(B_R)\le\rho_{d,\kappa}^2+u\}}
\,\mathrm du.
\]
All the integrands are nonnegative, so Tonelli's theorem
allows us to interchange the sum and the integral:
\[
R_{B_R}(\rho_{d,\kappa}^2+s^2)
=\int_0^{s^2}
\sum_{j=1}^{\infty}
\mathbf{1}_{\{\lambda_j(B_R)\le\rho_{d,\kappa}^2+u\}}
\,\mathrm du.
\]
The sum inside the integral is precisely the counting function.
Therefore,
\[
R_{B_R}(\rho_{d,\kappa}^2+s^2)
=\int_0^{s^2}
N_{B_R}(\rho_{d,\kappa}^2+u)\,\mathrm du.
\]		
Integrating \eqref{eq:ball-count-refined} yields
\[
R_{B_R}(\rho_{d,\kappa}^2+s^2)
\le L^{\mathrm{cl}}_{1,d}\operatorname{Vol}_\kappa(B_R)s^{d+2}
-\int_0^{s^2}\mathcal D_{d,\kappa}(R,\sqrt u)\,\mathrm du.
\]
The last integral is strictly positive, proving
\eqref{eq:ball-bly-riesz}.
For the sum bound, we sum \eqref{eq:ball-eigen} over $1\le j\le m$
and use
\[
\sum_{j=1}^m j^{2/d}
>\int_0^m x^{2/d}\,\mathrm dx
=\frac{d}{d+2}m^{1+2/d}.
\]
Weyl's law for a fixed ball proves sharpness of the leading
coefficient. The limit $\lambda_1(B_R)\to\rho_{d,\kappa}^2$ as
$R\to\infty$, proved above, shows that the shift cannot be increased.
\end{proof}

We next give three further consequences of the ball estimate.

\begin{corollary}\label{cor:sectors}
	Let $S_{R,\theta}\subset\mathbb{H}^2_\kappa$ {be a geodesic sector of radius $R>0$ and opening angle $0<\theta\le2\pi$.} We impose Dirichlet conditions on the outer arc and on the two radial sides. At $\theta=2\pi$, the sector is the disk with a radial slit. Set $c_\kappa=\kappa^2/4$ and
	\[
	q_\theta=\left\lceil \pi/\theta\right\rceil.
	\]
	Then, for every $s>0$,
	\begin{equation}\label{eq:sector-count-refined}
		N_{S_{R,\theta}}\left(c_\kappa+s^2\right)
		<\frac{{\operatorname{Vol}_\kappa}\left(S_{R,\theta}\right)}{4\pi}s^2
		-\frac1{2q_\theta}\left\lfloor \frac{sR}{\pi}+\frac14\right\rfloor.
	\end{equation}
	It also follows that
	\begin{equation}\label{eq:sector-count-linear}
		N_{S_{R,\theta}}\left(c_\kappa+s^2\right)
		<\frac{{\operatorname{Vol}_\kappa}\left(S_{R,\theta}\right)}{4\pi}s^2
		-\frac{R}{2\pi q_\theta}s+\frac3{8q_\theta}.
	\end{equation}
	In particular, the strict shifted P\'olya inequality holds, and
	\[
	\lambda_k\left(S_{R,\theta}\right)>c_\kappa+
	\frac{4\pi k}{{\operatorname{Vol}_\kappa}\left(S_{R,\theta}\right)},\qquad k\in\mathbb{N}.
	\]
\end{corollary}
\begin{proof}
It suffices to prove the count for $\kappa=1$. By separation of variables, the angular problem is
\[
-\Phi''(\varphi)=\nu^2\Phi(\varphi),
\qquad 0<\varphi<\theta,
\qquad \Phi(0)=\Phi(\theta)=0.
\]
Thus the angular eigenfunctions are
\[
\Phi_m(\varphi)=\sin\left(\frac{m\pi\varphi}{\theta}\right),
\qquad m\in\mathbb N.
\]
We substitute the angular eigenvalues \((m\pi/\theta)^2\) into the radial equation and get
\[
-w''+\frac{(m\pi/\theta)^2-1/4}{\sinh^2r}w=s^2w.
\]
Hence the corresponding radial orders are $\nu=m\pi/\theta$.

We first consider opening angles in the range $\pi\le\theta\le2\pi$,
for which $q_\theta=1$. Let $a_\theta=\pi/\theta\le1$. The simpler bound \eqref{eq:phasebound-basic} in Proposition~\ref{prop:phase} gives
\[
N_{S_{R,\theta}}\left(1/4+s^2\right)
\le\sum_{m\ge1}\left\lfloor g_{R,s}\left(a_\theta m\right)+\frac14\right\rfloor.
\]
Since each angular Dirichlet eigenvalue is simple, no multiplicity factor is needed here. The function $h\left(x\right)=g_{R,s}\left(a_\theta x\right)$ is nonnegative, nonincreasing, convex, and $1/2$-Lipschitz. By Theorem~\ref{thm:lattice-two}, we have
\[
\sum_{m\ge1}\left\lfloor h\left(m\right)+\frac14\right\rfloor
<\int_0^\infty h\left(x\right)\,\mathrm{d} x
-\frac12\left\lfloor h\left(0\right)+\frac14\right\rfloor.
\]
Using $h\left(0\right)=sR/\pi$ and substituting $\nu =a_\theta x$, we get
\[
N_{S_{R,\theta}}\left(1/4+s^2\right)
<\frac{\theta}{\pi}\int_0^\infty g_{R,s}\left(\nu \right)\,\mathrm{d} \nu 
-\frac12\left\lfloor \frac{sR}{\pi}+\frac14\right\rfloor.
\]
Applying \eqref{eq:volume-integral} with $d=2$ yields
\[
\frac{\theta}{\pi}\int_0^\infty g_{R,s}\left(\nu \right)\,\mathrm{d} \nu 
=\frac{\theta s^2}{4\pi}\int_0^R\sinh r\,\mathrm{d} r=\frac{s^2}{4\pi}{\operatorname{Vol}_1}\left(S_{R,\theta}\right).
\]
This proves \eqref{eq:sector-count-refined} when $\pi\le\theta\le2\pi$.

If $0<\theta<\pi$, let $q=q_\theta$. Then $\pi\le q\theta<2\pi$. The larger sector $S_{R,q\theta}$ is partitioned into $q$ rotated copies of $S_{R,\theta}$. Dirichlet bracketing \eqref{eq:bracketing} gives
\[
qN_{S_{R,\theta}}\left(\frac{1}{4}+s^2\right)
\le N_{S_{R,q\theta}}\left(\frac{1}{4}+s^2\right).
\]
Applying the estimate already proved to the larger sector, {we obtain}
\[
N_{S_{R,q\theta}}\left(\frac{1}{4}+s^2\right)
<\frac{q{\operatorname{Vol}_1}\left(S_{R,\theta}\right)}{4\pi}s^2
-\frac12\left\lfloor \frac{sR}{\pi}+\frac14\right\rfloor.
\]
Division by $q$ proves \eqref{eq:sector-count-refined}. For general $\kappa$, we use curvature scaling. Thus the same correction holds for every $\kappa>0$.

The elementary inequality
\[
\left\lfloor \frac{sR}{\pi}+\frac14\right\rfloor\ge\frac{sR}{\pi}-\frac34
\]
gives \eqref{eq:sector-count-linear}. Since the integer part in \eqref{eq:sector-count-refined} is nonnegative, we discard it to recover the strict shifted {P\'olya} count. Evaluating at an eigenvalue yields the stated eigenvalue inequality.
\end{proof}

\begin{corollary}\label{cor:unions}
Let $B_{R_j}\left(x_j\right), j\in\mathcal J$, be a nonempty finite or countable collection of disjoint geodesic balls in $\mathbb{H}^d_\kappa$. Suppose that
\[
V:=\sum_{j\in\mathcal J}\operatorname{Vol}_\kappa\left(B_{R_j}\right)<\infty.
\]
Define the open set
\[
\Omega=\bigsqcup_{j\in\mathcal J}B_{R_j}\left(x_j\right).
\]
Here $\bigsqcup$ denotes disjoint union. {Then the Dirichlet spectrum of $\Omega$ is discrete, and for every $s>0$,}
\[
N_\Omega\left(\rho_{d,\kappa}^2+s^2\right)<C_dVs^d.
\]
{Equivalently, for every $k\in\mathbb N$,}
\[
\lambda_k\left(\Omega\right)>\rho_{d,\kappa}^2+
4\pi^2\left(\frac{k}{\omega_dV}\right)^{2/d}.
\]	
\end{corollary}
\begin{proof}
The Dirichlet operator is the orthogonal direct sum of the ball operators \cite{Davies}. In the countably infinite case, finite total volume forces the component volumes to tend to zero. When $k=1$, Theorem~\ref{thm:balls} forces their first eigenvalues to tend to infinity. {{Thus only finitely many components contribute below any fixed spectral threshold. The counting function of $\Omega$ is}}
\[
{N_\Omega\left(\lambda\right)=\sum_{j\in\mathcal J}N_{B_{R_j}\left(x_j\right)}\left(\lambda\right).}
\]
Summing Theorem~\ref{thm:balls} gives the desired bound. Evaluating at $\lambda=\lambda_k\left(\Omega\right)$ proves the eigenvalue inequality.	
\end{proof}

\begin{corollary}\label{cor:tiles}
	Let $D\subset\mathbb{H}^d_\kappa$ be a nonempty open set.
	Suppose that a geodesic ball $B_R\subset\mathbb{H}^d_\kappa$ with $R>0$ contains $q\in\mathbb{N}$ pairwise disjoint open sets $D_1,\dots,D_q$, each isometric to $D$, such that
	\[
	\operatorname{Vol}_\kappa\left(B_R\setminus\bigcup_{j=1}^{q}D_j\right)=0.
	\]
	The copies of $D$ partition $B_R$ up to a set of zero hyperbolic volume.
	We impose Dirichlet boundary conditions on the entire boundary of $D$. Then, for every $s>0$,
	\[
	N_D\left(\rho_{d,\kappa}^2+s^2\right)
	<C_d\operatorname{Vol}_\kappa(D)s^d.
	\]
	Equivalently, for every $k\in\mathbb{N}$,
	\[
	\lambda_k(D)>
	\rho_{d,\kappa}^2+
	4\pi^2\left(
	\frac{k}{\omega_d\operatorname{Vol}_\kappa(D)}
	\right)^{2/d}.
	\]
\end{corollary}
\begin{proof}
Equation \eqref{eq:bracketing} implies \(qN_D\left(\lambda\right)\le N_{B_R}\left(\lambda\right)\). We have \(\operatorname{Vol}_\kappa\left(B_R\right)=q\operatorname{Vol}_\kappa\left(D\right)\). Applying Theorem~\ref{thm:balls} and dividing by $q$, we obtain the desired bound. This completes the proof.
\end{proof}

\begin{remark}
Corollary~\ref{cor:tiles} includes half balls, equal orthant sectors and radial chambers associated with finite reflection groups on $\mathbb{S}^{d-1}$. Corollary~\ref{cor:sectors} is stronger in dimension two because it imposes no finite tiling condition on the original sector.	
\end{remark}

\section[Eigenvalue sums on open sets of finite volume]{{Eigenvalue sums on open sets of finite volume}}\label{sec:general}

Throughout this section $d\ge2$ and $\Omega\subset\mathbb H^d_\kappa$
is a nonempty open set of finite volume
$V=\operatorname{Vol}_\kappa(\Omega)$. Recall the shift
$a_{d,\kappa}$ defined in \eqref{eq:averaged-shift}, and let
\[
\delta_{d,\kappa}:=\rho_{d,\kappa}^2-a_{d,\kappa}
=\begin{cases}
0,&d=2,\\[2pt]
(d-1)(d-3)\kappa^2/12,&d\ge3.
\end{cases}
\]

\subsection{Proof of Theorem~\ref{thm:averaged}}

For $\lambda>\rho_{d,\kappa}^2$, let
$t=\sqrt{\lambda-\rho_{d,\kappa}^2}$ and define
\[
\eta_{d,\kappa}(\lambda)=\frac{w_{d,\kappa}(t)}{2t},
\]
where $w_{d,\kappa}$ is the Plancherel density given in
\eqref{eq:plancherel-density}. We first establish an estimate for $\eta_{d,\kappa}(\lambda)$ that will be used below and {then prove Theorem~\ref{thm:averaged}.}

\begin{lemma}\label{lem:density-comparison}
For every $d\ge2$ and $\lambda>\rho_{d,\kappa}^2$,
\begin{equation}\label{eq:density-comparison}
\eta_{d,\kappa}(\lambda)
\le\frac d2 C_d(\lambda-a_{d,\kappa})^{d/2-1}.
\end{equation}
In dimension three this is an equality.
\end{lemma}
\begin{proof}
{First, note that}
\[
\frac d2 C_d=\frac{1}{(4\pi)^{d/2}\Gamma(d/2)}.
\]
If $d=2$, formula
\eqref{eq:planar-plancherel-density} gives
\[
\eta_{2,\kappa}(\lambda)
=\frac{1}{4\pi}\tanh\left(
\frac{\pi\sqrt{\lambda-\kappa^2/4}}{\kappa}\right)
\le\frac{1}{4\pi}.
\]
Since $C_2=1/(4\pi)$ and $d/2-1=0$, this is precisely
\eqref{eq:density-comparison} with $a_{2,\kappa}=\kappa^2/4$.

{For $d\ge3$, we consider odd and even dimensions separately.} Suppose first that $d=2q+1$, where $q\ge1$ and $q\in\mathbb N$. By \eqref{eq:pre-gamma-recurrence}, we have
\[
\left|\frac{\Gamma(q+\mathrm it/\kappa)}
{\Gamma(\mathrm it/\kappa)}\right|^2
=\frac{t^2}{\kappa^2}
\prod_{j=1}^{q-1}\left(\frac{t^2}{\kappa^2}+j^2\right).
\]
Consequently,
\begin{equation}\label{eq:odd-density}
\eta_{d,\kappa}(\lambda)
=\frac d2 C_d\,t
\prod_{j=1}^{q-1}(t^2+j^2\kappa^2).
\end{equation}
Let $r=q-1/2=d/2-1$.
Applying the weighted arithmetic--geometric mean inequality to
$t^2$ with weight $1/(2r)$ and to each
$t^2+j^2\kappa^2$, $1\le j\le q-1$, with weight $1/r$, we obtain
\[
t\prod_{j=1}^{q-1}(t^2+j^2\kappa^2)
\le\left(t^2+\frac{\kappa^2}{r}\sum_{j=1}^{q-1}j^2\right)^r.
\]
Since
\[
\frac1r\sum_{j=1}^{q-1}j^2
=\frac{q(q-1)}3
=\frac{(d-1)(d-3)}{12},
\]
the expression in parentheses is
$t^2+\delta_{d,\kappa}=\lambda-a_{d,\kappa}$.
This proves \eqref{eq:density-comparison} in odd dimensions.
When $q=1$, the single weighted factor gives equality.

Now suppose that $d=2q+2$, where $q\ge1$ and $q\in\mathbb N$. By \eqref{eq:pre-gamma-recurrence} and \eqref{eq:pre-gamma-imaginary-modulus}, we obtain
\[
\eta_{d,\kappa}(\lambda)
=\frac d2 C_d\tanh(\pi t/\kappa)
\prod_{j=1}^{q}\left(t^2+(j-1/2)^2\kappa^2\right).
\]
Here $0<\tanh(\pi t/\kappa)<1$. The arithmetic--geometric mean inequality gives
\[
\prod_{j=1}^{q}\left(t^2+(j-1/2)^2\kappa^2\right)
\le\left(t^2+\frac{\kappa^2}{q}
\sum_{j=1}^{q}(j-1/2)^2\right)^q.
\]
Since $q=d/2-1$ and
\[
\frac1q\sum_{j=1}^{q}(j-1/2)^2
=\frac{(d-1)(d-3)}{12},
\]
{this is exactly the desired comparison.}
\end{proof}

\begin{proof}[Proof of Theorem~\ref{thm:averaged}]
If $\lambda\le\rho_{d,\kappa}^2$, then $R_\Omega(\lambda)=0$
and \eqref{eq:general-riesz} follows immediately.
{When $d=2$ and $\lambda>\kappa^2/4$, set $s=\sqrt{\lambda-\kappa^2/4}>0$.} Lemma~\ref{lem:riesz} and
\eqref{eq:planar-plancherel-density} give
\[
R_\Omega(\lambda)
\le\frac{V}{2\pi}{\int_0^{s}(s^2-t^2)t\tanh(\pi t/\kappa)\,\mathrm dt.}
\]
Using $0<\tanh(\pi t/\kappa)<1$ and
\[
{\int_0^{s}(s^2-t^2)t\,\mathrm dt=\frac{s^4}{4},}
\]
we obtain \eqref{eq:planar-riesz}, which is
\eqref{eq:general-riesz} for $d=2$ because
$L^{\mathrm{cl}}_{1,2}=1/(8\pi)$.

For $d\ge3$ and $\lambda>\rho_{d,\kappa}^2$, by the change of variables $\xi=\rho_{d,\kappa}^2+t^2$ in \eqref{eq:riesz-estimate}, we obtain
\[
R_\Omega(\lambda)
\le V\int_{\rho_{d,\kappa}^2}^{\lambda}
(\lambda-\xi)\eta_{d,\kappa}(\xi)\,\mathrm d\xi.
\]
By Lemma~\ref{lem:density-comparison} and
$a_{d,\kappa}\le\rho_{d,\kappa}^2$, we obtain
\[
R_\Omega(\lambda)
\le\frac d2 C_d V
\int_{a_{d,\kappa}}^{\lambda}
(\lambda-\xi)(\xi-a_{d,\kappa})^{d/2-1}\,\mathrm d\xi.
\]
Set $z=\lambda-a_{d,\kappa}$ and $v=\xi-a_{d,\kappa}$.
Direct integration gives
\[
\int_0^z(z-v)v^{d/2-1}\,\mathrm dv
=\frac{z^{1+d/2}}{(d/2)(d/2+1)}.
\]
Multiplication by $(d/2)C_dV$ proves
\eqref{eq:general-riesz}, because
$C_d/(d/2+1)=L^{\mathrm{cl}}_{1,d}$.

For every $\lambda\in\mathbb R$ and integer $m\ge1$, it is easy to see that
\[
\sum_{j=1}^m\lambda_j(\Omega)\ge m\lambda-R_\Omega(\lambda).
\]
Taking $\lambda=a_{d,\kappa}+z$, $z\ge0$, yields
\[
\sum_{j=1}^m\lambda_j(\Omega)
\ge a_{d,\kappa}m+mz-L^{\mathrm{cl}}_{1,d}Vz^{1+d/2}.
\]
{The sum of the last two terms, viewed as a function of $z$, attains its maximum at}
\[
z_{*}=\left(\frac{m}{C_dV}\right)^{2/d}
=4\pi^2\left(\frac{m}{\omega_dV}\right)^{2/d},
\]
since $(1+d/2)L^{\mathrm{cl}}_{1,d}=C_d$. At this value, the critical point equation gives
\[
m=(1+d/2)L^{\mathrm{cl}}_{1,d}Vz_{*}^{d/2}.
\]
Hence,
\[
mz_{*}-L^{\mathrm{cl}}_{1,d}Vz_{*}^{1+d/2}
=\frac{d}{d+2}mz_{*}.
\]
Substituting this value into the preceding lower bound yields
\[
\sum_{j=1}^m\lambda_j(\Omega)
\ge a_{d,\kappa}m+
\frac{d}{d+2}\,4\pi^2
\left(\frac{m}{\omega_dV}\right)^{2/d}m,
\]
which proves \eqref{eq:sum-inequality} in every dimension $d\ge2$. For $d=2$, the maximizing value is $z_{*}=4\pi m/V$, and the
resulting bound is exactly \eqref{eq:planar-sum}. Conversely, since
\[
R_\Omega(\lambda)
=\sup_{m\in\mathbb N_0}
\left(m\lambda-\sum_{j=1}^m\lambda_j(\Omega)\right),
\]
substituting \eqref{eq:sum-inequality} and enlarging the supremum to all real
$m\ge0$ gives \eqref{eq:general-riesz} by the same optimization.
This proves the asserted equivalence.

For a fixed smooth bounded open set, Weyl's law gives
\[
\sum_{j=1}^m\lambda_j(\Omega)
\sim\frac{d}{d+2}\,4\pi^2
(\omega_dV)^{-2/d}m^{1+2/d}\qquad m\to\infty.
\]
The linear term is of smaller order. Thus the leading coefficient
in \eqref{eq:sum-inequality} cannot be increased uniformly. To prove optimality of the planar shift, suppose that a constant
$b>\kappa^2/4$ could replace $\kappa^2/4$ in
\eqref{eq:planar-sum} {for every open set of finite area and every $m\in\mathbb N$.}
Taking $m=1$ and $\Omega=B_R\subset\mathbb H^2_\kappa$ would give
\[
\lambda_1(B_R)\ge b+\frac{2\pi}{\operatorname{Vol}_\kappa(B_R)}.
\]
As $R\to\infty$, the volume tends to infinity and
$\lambda_1(B_R)\to\kappa^2/4$, as proved in
Section~\ref{sec:balls}. {This contradicts $b>\kappa^2/4$.}
\end{proof}

\subsection{The global estimate proposed by Strichartz}\label{subsec:strichartz}

In this subsection we assume $d\ge3$. Recall that the spectral variables $\lambda$ and $t$ are related by
$\lambda=\rho_{d,\kappa}^2+t^2$, with $t\ge0$.
We now define the renormalized eigenvalue sum used in
\cite[Section~5]{Strichartz}, with the normalization of this paper.
For $E\ge\rho_{d,\kappa}^2$, let
\begin{equation}\label{eq:ambient-count}
\mathcal A_{d,\kappa}(E)
=\int_{\rho_{d,\kappa}^2}^{E}\eta_{d,\kappa}(\lambda)\,\mathrm d\lambda.
\end{equation}
$\eta_{d,\kappa}(\lambda)$ is positive for $\lambda>\rho_{d,\kappa}^2$, so this function
is continuous and strictly increasing from zero to infinity. Therefore, we can define the renormalized eigenvalue sum as follows.
\begin{definition}
	\label{def:renormalized-sum}
	For each $x>0$, let $E_x>\rho_{d,\kappa}^2$ be the unique number
	satisfying $\mathcal A_{d,\kappa}(E_x)=x$. The renormalized eigenvalue sum is
	\begin{equation}\label{eq:renormalized-sum}
		\Sigma_{d,\kappa}(x)
		=\int_{\rho_{d,\kappa}^2}^{E_x}
		\lambda\eta_{d,\kappa}(\lambda)\,\mathrm d\lambda.
	\end{equation}
\end{definition}
By changing variables $y=\mathcal A_{d,\kappa}(\lambda)$, we obtain
\begin{equation}\label{eq:inverse-integral}
\Sigma_{d,\kappa}(x)=\int_0^x E_y\,\mathrm dy.
\end{equation}
For completeness, we compare this definition with
\cite[equations~(5.12)--(5.13)]{Strichartz}. First take $\kappa=1$, and set $\rho=(d-1)/2$. Let \(h_{\mathrm S}\), \(H_{\mathrm S}\), and \(\Sigma_{\mathrm S}\) denote the functions $h, H$, and \(\Sigma\) from that paper, and let \(\sigma_{d-1}=|\mathbb S^{d-1}|=d\omega_d\). Strichartz's $\omega_{d-1}$ denotes this sphere area, whereas our $\omega_{d-1}$ denotes the volume of the Euclidean unit ball in dimension $d-1$. His formulas become
\[
H_{\mathrm S}(r)=\int_0^r h_{\mathrm S}(t)\,\mathrm dt,
\qquad
\Sigma_{\mathrm S}(x)=\sigma_{d-1}
\int_0^{H_{\mathrm S}^{-1}(x/\sigma_{d-1})}
(\rho^2+t^2)h_{\mathrm S}(t)\,\mathrm dt.
\]
By the change of variables $\lambda=\rho^2+t^2$, we obtain
\[
\eta_{d,1}(\rho^2+t^2)\,\mathrm d\lambda
=\sigma_{d-1}h_{\mathrm S}(t)\mathrm dt,\qquad t>0.
\]
Consequently,
\[
\begin{aligned}
	\mathcal A_{d,1}(E)
	&=\int_{\rho^2}^{E}
	\eta_{d,1}(\lambda)\,\mathrm d\lambda\\
	&=\sigma_{d-1}\int_0^{\sqrt{E-\rho^2}}h_{\mathrm S}(t)\,\mathrm dt\\
	&=\sigma_{d-1}H_{\mathrm S}\bigl(\sqrt{E-\rho^2}\bigr).
\end{aligned}
\]
Thus, the relation $\mathcal A_{d,1}(E_x)=x$ implies $\sqrt{E_x-\rho^2}=H_{\mathrm S}^{-1}(x/\sigma_{d-1})$. Using the same change of variables, we obtain
\[
\begin{aligned}
	\Sigma_{d,1}(x)
	&=\int_{\rho^2}^{E_x}
	\lambda\,\eta_{d,1}(\lambda)\,\mathrm d\lambda\\
	&=\sigma_{d-1}\int_0^{H_{\mathrm S}^{-1}(x/\sigma_{d-1})}
	(\rho^2+t^2)h_{\mathrm S}(t)\,\mathrm dt
	=\Sigma_{\mathrm S}(x).
\end{aligned}
\]
This proves that the two definitions coincide when $\kappa=1$.
For general \(\kappa>0\), metric scaling yields the corresponding relations.
\begin{remark}
Our notation for the Weyl constant differs from that of
Strichartz \cite{Strichartz}. We write $C_d=\omega_d(2\pi)^{-d}$ for the coefficient in the Weyl asymptotic formula for the
counting function. The constant denoted by $C_d$ in
\cite{Strichartz} corresponds instead to $C_d^{-2/d}=4\pi^2\omega_d^{-2/d}$ in our notation.	
\end{remark}

\begin{proposition}\label{prop:strichartz}
For every integer $d\ge3$, every $\kappa>0$, and every $x>0$,
\begin{equation}\label{eq:strichartz-global}
a_{d,\kappa}x
\le\Sigma_{d,\kappa}(x)
-\frac{d}{d+2}\,4\pi^2\omega_d^{-2/d}x^{1+2/d}
\le\rho_{d,\kappa}^2x.
\end{equation}
The lower coefficient $a_{d,\kappa}$ cannot be increased, and the
upper coefficient $\rho_{d,\kappa}^2$ cannot be decreased.
When $d=3$, both inequalities are equalities for every $x>0$.
For $\kappa=1$, this proves the estimate proposed in
\cite[equation~(5.24)]{Strichartz}, with its conjectured constants.
\end{proposition}
\begin{proof}
We first prove the lower estimate. Define, for every real $E$,
\[
\mathcal R_{d,\kappa}(E)
=\int_{\rho_{d,\kappa}^2}^{\infty}
(E-\xi)_+\eta_{d,\kappa}(\xi)\,\mathrm d\xi.
\]
For $E>\rho_{d,\kappa}^2$, differentiation under the integral gives
$\mathcal R'_{d,\kappa}(E)=\mathcal A_{d,\kappa}(E)$.
Thus the function $xE-\mathcal R_{d,\kappa}(E)$ increases when $E<E_x$ and
decreases when $E>E_x$. It attains its maximum at $E_x$. For $E\le\rho_{d,\kappa}^2$ it equals $xE$. Substituting $E=E_x$, we obtain
\[
\begin{aligned}
	\sup_{E\in\mathbb R}\{xE-\mathcal R_{d,\kappa}(E)\}
	&=xE_x-\int_{\rho_{d,\kappa}^2}^{E_x}
	(E_x-\lambda)\eta_{d,\kappa}(\lambda)\,\mathrm d\lambda\\
	&=E_x\left(
	x-\int_{\rho_{d,\kappa}^2}^{E_x}
	\eta_{d,\kappa}(\lambda)\,\mathrm d\lambda
	\right)
	+\int_{\rho_{d,\kappa}^2}^{E_x}
	\lambda\eta_{d,\kappa}(\lambda)\,\mathrm d\lambda.
\end{aligned}
\]
Since $\mathcal A_{d,\kappa}(E_x)=x$, we obtain
\begin{equation}\label{eq:ambient-legendre}
\sup_{E\in\mathbb R}\{xE-\mathcal R_{d,\kappa}(E)\}
=\Sigma_{d,\kappa}(x).
\end{equation}
Lemma~\ref{lem:density-comparison}, as in the proof of Theorem~\ref{thm:averaged}, yields
\[
\mathcal R_{d,\kappa}(E)
\le L^{\mathrm{cl}}_{1,d}(E-a_{d,\kappa})_+^{1+d/2}.
\]
We next choose $E=a_{d,\kappa}+(x/C_d)^{2/d}$ in
\eqref{eq:ambient-legendre}. The same optimization as before gives
\[
\Sigma_{d,\kappa}(x)
\ge a_{d,\kappa}x
+\frac{d}{d+2}C_d^{-2/d}x^{1+2/d},
\]
which is the first inequality in \eqref{eq:strichartz-global}.

For the upper estimate, we claim
\begin{equation}\label{eq:density-lower}
\eta_{d,\kappa}(\lambda)
\ge\frac d2 C_d(\lambda-\rho_{d,\kappa}^2)^{d/2-1},
\qquad \lambda>\rho_{d,\kappa}^2.
\end{equation}
If $d=2q+1$, where $q\ge1$ and $q\in\mathbb N$, formula~\eqref{eq:odd-density} proves this immediately:
each factor $t^2+j^2\kappa^2$ is at least $t^2$ and hence
\[
t\prod_{j=1}^{q-1}(t^2+j^2\kappa^2)\ge t^{2q-1}.
\]
If $d=2q+2$, where $q\ge1$ and $q\in\mathbb N$, set $z=t/\kappa>0$.
The elementary inequality
\[
e^{2\pi z}>1+\frac{(2\pi z)^2}{2}>1+8z^2
\]
implies $(z^2+1/4)\tanh(\pi z)>z^2$. Indeed, substituting
$\tanh(\pi z)=(e^{2\pi z}-1)/(e^{2\pi z}+1)$ shows that the last
inequality is equivalent to $e^{2\pi z}>1+8z^2$. Thus, we obtain
\[
\tanh(\pi t/\kappa)
\prod_{j=1}^{q}\bigl(t^2+(j-1/2)^2\kappa^2\bigr)
\ge t^{2q}.
\]
The even-dimensional density formula in the proof of
Lemma~\ref{lem:density-comparison} now proves \eqref{eq:density-lower}. Integrating \eqref{eq:density-lower} yields
\[
\mathcal A_{d,\kappa}(E)
\ge C_d(E-\rho_{d,\kappa}^2)^{d/2},
\qquad E\ge\rho_{d,\kappa}^2.
\]
At $E=E_y$, this gives
\[
E_y\le\rho_{d,\kappa}^2+(y/C_d)^{2/d}.
\]
Integrating over $0<y<x$ and using \eqref{eq:inverse-integral}
proves the upper estimate in \eqref{eq:strichartz-global}.

We next verify optimality of the two linear coefficients.
As $x\rightarrow0$, by the definition of $E_x$, we have
$E_x\to\rho_{d,\kappa}^2$. Integrating over $[\rho_{d,\kappa}^2,E_x]$, we obtain
\[
\rho_{d,\kappa}^2
\int_{\rho_{d,\kappa}^2}^{E_x}
\eta_{d,\kappa}(\lambda)\,\mathrm{d}\lambda
\le
\int_{\rho_{d,\kappa}^2}^{E_x}
\lambda\eta_{d,\kappa}(\lambda)\,\mathrm{d}\lambda
\le
E_x\int_{\rho_{d,\kappa}^2}^{E_x}
\eta_{d,\kappa}(\lambda)\,\mathrm{d}\lambda.
\]
Using the definitions of $\Sigma_{d,\kappa}(x)$ and $E_x$, we obtain
\[
\rho_{d,\kappa}^2
\le\frac{\Sigma_{d,\kappa}(x)}{x}\le E_x.
\]
Consequently,
\begin{equation}\label{eq:strichartz-small-limit}
\lim_{x\rightarrow0}
\frac{\Sigma_{d,\kappa}(x)
-\frac{d}{d+2}C_d^{-2/d}x^{1+2/d}}{x}
=\rho_{d,\kappa}^2.
\end{equation}
This proves optimality of the upper coefficient. By {Strichartz's} asymptotic formula \cite[equation~(5.16)]{Strichartz} and the scaling relation \cite[equation~(5.25)]{Strichartz}, we obtain
\[
\lim_{x\to\infty}
\frac{\Sigma_{d,\kappa}(x)
	-\frac{d}{d+2}C_d^{-2/d}x^{1+2/d}}{x}
=a_{d,\kappa}.
\]
Hence the lower coefficient cannot be increased. Finally, when $d=3$, $a_{3,\kappa}=\rho_{3,\kappa}^2=\kappa^2$,
so the two estimates coincide for every $x>0$.
\end{proof}

\begin{remark}
{The lower bound of \cite{Strichartz}, which also follows from Lemma~\ref{lem:riesz} and \eqref{eq:ambient-legendre}, is}
\[
\sum_{j=1}^m\lambda_j(\Omega)
\ge V\Sigma_{d,\kappa}\left(\frac{m}{V}\right).
\]
{However, the upper bound in Proposition~\ref{prop:strichartz} gives}
$$
m\rho_{d,\kappa}^2 + \frac{dm}{d+2}\,4\pi^2\omega_d^{-2/d}\left(\frac{m}{V}\right)^{2/d}
\ge V\Sigma_{d,\kappa}\left(\frac{m}{V}\right).
$$
{We cannot deduce the desired conclusion of Conjecture ~\ref{conj:laptev-bly} from these two relations. Thus the full shift in
Conjecture~\ref{conj:laptev-bly} for general domains in dimensions
$d\ge4$ is not established here.}	
\end{remark}

\section*{Acknowledgments}
This work was supported by the National Natural Science Foundation of China (Grant No.~12371110).

\section*{Statements and Declarations}
\medskip
\noindent\textbf{Author Contributions}
{The authors contributed equally to this work.}

\medskip
\noindent\textbf{Competing interests}
There are no competing interests.

\medskip
\noindent\textbf{Data availability}
Data sharing is not applicable to this article as no datasets
were generated or analyzed during the current study.

\medskip
\noindent\textbf{Ethical approval}
{Not applicable.
}

\medskip
\noindent\textbf{Declaration of generative AI}
AI tools (GPT-6 Astra) were used to assist with language presentation and
LaTeX formatting. The authors take full responsibility for
the mathematical content and the final manuscript.

\end{document}